\documentclass[11pt]{article}
\usepackage[margin=1in]{geometry}
\usepackage{amsmath,amssymb,amsthm,mathtools}
\usepackage[T1]{fontenc}
\usepackage{lmodern}
\usepackage{microtype}
\usepackage[hidelinks]{hyperref}
\newtheorem{theorem}{Theorem}
\newtheorem{lemma}[theorem]{Lemma}
\newtheorem{proposition}[theorem]{Proposition}

\theoremstyle{remark}

\theoremstyle{plain}
\newcommand{\CP}{\mathbb{CP}}
\newcommand{\R}{\mathbb R}
\newcommand{\C}{\mathbb C}
\newcommand{\Z}{\mathbb Z}
\newcommand{\Q}{\mathbb Q}
\newcommand{\F}{\mathbb F}
\newcommand{\Tr}{\operatorname{Tr}}
\newcommand{\Spec}{\operatorname{Spec}}
\newcommand{\Vol}{\operatorname{Vol}}
\newcommand{\Ric}{\operatorname{Ric}}

\title{A proof of the spectral rigidity for complex projective space}
\author{Kefeng Liu \\ Chongqing University of Technology\\ kefengliu@cqut.edu.cn
\and 
Hao Xu \\ Chongqing University of Technology\\ haoxu@cqut.edu.cn}
\date{}
\begin{document}
\maketitle
\begin{abstract}
We give a proof of a folklore conjecture that a compact 
K\"ahler manifold with the same scalar Laplace spectrum as complex
projective space equipped with the Fubini--Study metric is
holomorphically isometric to it. More generally, we prove that the
spectrum of the Hodge Laplacian in any one fixed degree $p$ determines
$\CP^n$ up to holomorphic isometry among compact connected K\"ahler
manifolds, for every fixed pair of integers $(p,n)$ with $0\leq p\leq n$.
Equality of a single Hodge spectrum first forces periodic geodesic flow
and the integral cohomology ring of complex projective space.
The ordinary first heat coefficient then identifies the total scalar
curvature except on a Pell sequence of pairs. For those pairs, a clean
wave-trace argument, a Krawtchouk estimate, and a mod-two
projective-bundle obstruction determine the geodesic return holonomy.
A canonical K\"ahler involution recovers a weighted first heat coefficient
which is nonzero throughout the exceptional sequence. The
Kobayashi--Ochiai theorem and a fixed-class second-heat-coefficient
calculation complete the argument.
\end{abstract}

\section{Introduction}

We use the nonnegative Hodge Laplacian
$\Delta_p=d\delta+\delta d$ on real $p$-forms, complexifying its
eigenspaces when useful. Complexification preserves their dimensions
and hence the spectral multiplicities. The case $n=0$ is the one-point
case. Assume henceforth $n\geq1$. Weyl's law determines the real
dimension from the spectrum, so $\dim_\R M=2n$.

Rescaling both metrics by the same positive constant rescales both
spectra by its reciprocal and preserves the hypothesis. We therefore
normalize the model to have holomorphic sectional curvature $4$.
With the ordinary metric K\"ahler form
$\omega(X,Y)=g(JX,Y)$ and the positive hyperplane generator
$h\in H^2(\CP^n;\Z)$, the model has
\[
 [\omega_0]=\pi h,\qquad
 V_0=\frac{\pi^n}{n!},\qquad s_0=4n(n+1).
\]
Throughout, $s$ denotes the Riemannian scalar curvature and
$dV=\omega^n/n!$. Put
\begin{equation}\label{notation}
 D=2n,\qquad N=\binom Dp,\qquad
 q=p^2-2np+\frac{n(2n-1)}3.
\end{equation}
Our heat-coefficient convention is
\begin{equation}\label{heatconvention}
 \Tr(e^{-t\Delta_p})\sim(4\pi t)^{-D/2}
       \sum_{j\geq0}a_{j,p}t^j\qquad(t\downarrow0).
\end{equation}
In particular $a_{0,p}=N\Vol(M)$, so isospectrality also fixes volume.

\begin{theorem}\label{main}
For fixed $(p,n)$ with $0\leq p\leq n$, where $p,n$ are integers,
let $(M,g,J)$ be a compact connected K\"ahler manifold without boundary
such that
\[
 \Spec_p(M,g)=\Spec_p(\CP^n,g_0),
\]
including multiplicities, where $g_0$ is a Fubini--Study metric.
Then $(M,g,J)$ is holomorphically isometric to $(\CP^n,g_0,J_0)$.
\end{theorem}

The hypothesis concerns the spectrum in the \emph{one fixed degree $p$};
it does not require simultaneous isospectrality in several degrees.
The results from the literature are cited separately from the
additional argument supplied here.

The proof has four stages.
\begin{enumerate}
\item \emph{From spectral arithmetic to topology.}
The explicit Fubini--Study eigenvalue branches imply
$e^{2\pi i\sqrt\lambda}\to1$ at high energy. Hence
$U_p(2\pi)-I$ is compact. Egorov's theorem, tested against arbitrary
scalar principal symbols, forces every geodesic to return after time
$2\pi$. The universal cover is then compact. The integral
Bott--Samelson--McCleary theorem and the nonzero K\"ahler class give
projective-space cohomology; the Lefschetz fixed-point theorem excludes
nontrivial deck transformations. Thus $\pi_1(M)=0$ and
$H^*(M;\Z)=\Z[x]/(x^{n+1})$. Equal volume fixes $[\omega]=\pi x$.

\item \emph{From total scalar curvature to the complex manifold.}
Write $c_1(M)=kx$. The scalar-curvature pairing gives
$V^{-1}\int_Ms\,dV=4nk$. When $q\ne0$, the coefficient $a_{1,p}$
has a nonzero scalar-curvature multiplier, so isospectrality yields
$k=n+1$. The integral K\"ahler class $x$ is ample; the
Kobayashi--Ochiai maximal-index theorem identifies $M$ biholomorphically
with $\CP^n$, with the model K\"ahler class.

\item \emph{Recovering the missing invariant on the Pell sequence.}
When $q=0$, the ordinary $a_{1,p}$ vanishes identically. A nonzero
leading model wave-trace coefficient first forces the shorter return
$G^\pi=I$. The exterior holonomy character, a Krawtchouk estimate,
and a mod-two obstruction on $P_\R(TM)$ then determine the entire
return holonomy. It agrees, up to parity, with the principal symbol
of an exact K\"ahler involution $\mathcal K_p$ commuting with
$\Delta_p$. Compactness of $U_p(\pi)-\varepsilon\mathcal K_p$ and
separation of the two signs determine the weighted multiplicities
outside finitely many eigenvalues. The resulting weighted first heat
coefficient has a nonzero scalar-curvature multiplier and recovers
$k=n+1$ also in every exceptional case.

\item \emph{From the fixed complex manifold and class to the metric.}
The two Apte identities fix the relevant combinations of curvature
integrals in one K\"ahler class. Substitution in $a_{2,p}$ leaves the
positive multiple $\binom{2n}{p}/60$ of scalar-curvature variance.
Its vanishing forces constant scalar curvature. The same identities
and the orthogonal K\"ahler curvature decomposition then force both
the traceless Ricci and Bochner tensors to vanish. The metric has
holomorphic sectional curvature $4$, and complex-space-form
uniformization gives the required holomorphic isometry.
\end{enumerate}

The nonexceptional proof uses Sections~\ref{sec:analytic}--\ref{sec:topology}
and then Section~\ref{sec:metric}; only the first two subsections of
Section~\ref{sec:analytic} are needed in that route.
Sections~\ref{sec:return}--\ref{sec:recovery} supply the additional
argument for the Pell sequence. Section~\ref{sec:consequences}
records consequences of the same proof.

\paragraph{The scalar case.}
For $p=0$, one has $N=1$ and $\beta_0=1/6$, so there is no exceptional
dimension. The model eigenvalues are simply $4k(k+n)$, and the above
first, second, and fourth stages prove scalar spectral rigidity in
every complex dimension. In particular, the clean return trace,
holonomy-character estimate, projective-bundle obstruction, and
weighted involution are unnecessary for this case. The full function
spectrum is still essential to the first stage: its arithmetic
high-energy behavior supplies periodicity and topology, information
not provided by the first three ordinary heat coefficients alone.

\

\noindent{\bf Acknowledgements} This work was supported by the
National Key Research and Development Program of China (NKPs) [Grant
Number 2024YFA1013201], the Natural Science Foundation of Chongqing
(NSFCQ) [Grant Number CSTB2024NSCQ-LZX0040], the Special Project of
Chongqing Municipal Science and Technology Bureau [Grant Number
2025CCZ015]. 

\

\noindent{\bf Disclosure of AI assistance.}
GPT-6 Astra was used as an auxiliary tool for literature search, proofreading, 
and checks of mathematical arguments.

\section{Historical context and relation to earlier results}

The question belongs to the inverse spectral problem: how much of a
Riemannian metric is determined by the spectrum of its Laplacian?  Milnor's
isospectral nonisometric flat tori show that no unrestricted uniqueness
theorem is possible \cite{Milnor}.  For highly symmetric reference metrics,
however, the heat expansion provides strong rigidity conditions.
The work of McKean--Singer and Patodi identified the curvature information
encoded in the early heat coefficients, including the Hodge Laplacians on
forms \cite{McKeanSinger,Patodi,PatodiEigenforms}.  The monograph of
Berger--Gauduchon--Mazet is a classical account of the spectral setting
\cite{BergerGauduchonMazet}.  Gilkey's monograph gives
a systematic treatment of these local invariants \cite{Gilkey}.

For complex projective space, the assertion that the scalar Laplace
spectrum characterizes the Fubini--Study metric among compact connected
K\"ahler manifolds was a longstanding folklore conjecture. The general question for
the $p$-form spectrum was first stated explicitly by Chen--Vanhecke
\cite{ChenVanhecke} in 1980.

Tanno proved rigidity from the function
spectrum in complex dimensions $n\leq6$ and from the one-form spectrum
for $8\leq n\leq51$ \cite{Tanno1973,Tanno1974}.  The function case $n=6$
uses information beyond the first three heat coefficients; the direct
nonnegative-curvature-decomposition argument using $a_0,a_1,a_2$ covers
$n\leq5$ \cite[Remark 2.6]{HLXZ}.  These are statements about the ranges
established in the cited works, not an assertion that no other approach
can improve them.

Chen--Vanhecke established two-form rigidity except in complex dimensions
$2$ and $8$ \cite{ChenVanhecke}.  Goldberg treated the remaining dimensions,
and Gauchman--Goldberg and Perrone subsequently proposed related
extensions or alternate proofs \cite{Goldberg,GauchmanGoldberg,Perrone}.
The qualifications in this history matter.  Li showed that Goldberg's
dimension-two argument is valid, but identified a gap in the volume
argument for dimension eight, inherited by the later
Gauchman--Goldberg paper, and a scale-normalization problem in Perrone's
argument \cite[Introduction and Section 2.3]{Li}.  Thus those earlier
claims should not be counted as independent complete proofs in every
dimension.  Li supplied the missing rigorous conclusion, obtaining
two-form rigidity in all dimensions \cite[Corollary 1.4]{Li}.
Fujita's sharp volume bound for K\"ahler--Einstein Fano manifolds is a
central input in Li's treatment of the exceptional curvature-coefficient
case \cite{Fujita}.

Li \cite[Theorem 1.2]{Li} proved the rigidity
when $p$ is positive and even, $2\leq p\leq2(n-1)$, and
\[
 p^2-2np+\frac{n(2n-1)}3\ne0.
\]
The excluded integer pairs form a Pell sequence.  In the convention
$p\leq n$ used here, the positive even-degree exceptions begin with
$(p,n)=(20,48),(3976,9408),\ldots$.  They are exceptions to the theorem's
proof, not known counterexamples.  The additional endpoint $p=2n$ in
Li's broader formulation is equivalent to the scalar problem by Hodge
duality and is outside the positive-degree range $p\leq n$ of the present
question.

Li--Sun--Zhu's work on the spectral rigidity of Ricci soliton and
Einstein-type manifolds \cite{LiSunZhu} develops related characterizations under
additional geometric hypotheses. In the K\"ahler setting they assume
that the metric is cohomologically Einstein, meaning that $c_1(M)$ is
proportional to its K\"ahler class. Equality of the first three heat
coefficients with those of Fubini--Study space then gives holomorphic
isometry for $p=0$ and $p=2$ in all dimensions, and for $p\geq3$ outside
the exceptional Pell pairs.
Their paper incorporates earlier K\"ahler results of Li and also treats
gradient shrinking Ricci solitons in the Riemannian setting and suitable
almost-isospectral extensions. This work places the cohomologically
Einstein case in a broader spectral-rigidity framework and is closely
related to the scalar-curvature and Chern--Weil arguments used below.

The work of Huang--Liu--Xu--Zhi 
\cite[Theorem~2.10]{HLXZ} gives explicit dimension ranges
for every fixed $p\geq3$ by analyzing the signs of the first three heat
coefficients.  For example, their method gives $4\leq n\leq6$ and
$10\leq n\leq136$ for $p=3$.  They also explain why the zero set of the
first scalar-curvature coefficient lies outside those ranges
\cite[Remark 2.12 and Examples 2.13--2.14]{HLXZ}.  

The argument in the present note uses the full spectrum to
obtain additional geometric and weighted spectral information; it is
therefore not limited to the positivity conditions of the ordinary first
three heat coefficients.

The present work was inspired by Ping Li's work on spectral rigidity
\cite{Li}.
His use of the Hodge heat coefficients, K\"ahler curvature identities,
and the cohomologically Einstein condition guides the reductions here.
Our additional step is to use the full fixed-degree spectrum to obtain
periodicity and topology, and then to recover the missing weighted
coefficient in the exceptional cases.

The other inputs have separate origins.  Ikeda--Taniguchi computed the
Hodge spectra of spheres and complex projective spaces
\cite{IkedaTaniguchi}.  The relation between spectra and closed
bicharacteristics is developed by Duistermaat--Guillemin, with the
pseudodifferential and Fourier integral calculus treated systematically
by H\"ormander \cite{DG,HormanderIV}; Seeley's complex powers provide the
standard elliptic-functional-calculus background \cite{Seeley}.
For manifolds with all geodesics closed, Besse's monograph and McCleary's
loop-space results provide the cohomological restrictions used here
\cite{Besse,McCleary}.  The precise integral-cohomology conclusion for
simply connected Besse manifolds is stated in the introduction of
Radeschi--Wilking \cite{RadeschiWilking}; this is the formulation applied
in this note.  Wadsley's theorem is part of the general background on
common geodesic periods \cite{Wadsley}.

Finally, the passage from a positive integral K\"ahler class to an ample
line bundle uses Kodaira's theory \cite{Kodaira}, and the resulting
maximal-index characterization is due to Kobayashi--Ochiai
\cite{KobayashiOchiai}.  The final metric comparison uses Apte's
Chern--Weil identities \cite{Apte}.  Standard background for the K\"ahler
identities and Hodge theory is given in Huybrechts \cite{Huybrechts}; the
characteristic-class and projective-bundle calculations use the
conventions of Milnor--Stasheff \cite{MilnorStasheff}.  The final
classification of complete simply connected complex space forms is part
of the standard differential geometry.

\section{Detailed analytic preliminaries}
\label{sec:analytic}

Throughout this section $X$ is a closed connected Riemannian manifold of real
dimension $D\geq2$, and $E\to X$ is a Hermitian vector bundle. Operators and
forms may be complexified without changing the real multiplicities of the
Hodge spectrum. We take the nonnegative convention for the Laplacian.
The finite-dimensional zero eigenspace never affects a principal symbol or
any negative-power heat coefficient. Thus, when applying the calculus of
elliptic square roots, one can replace $\Delta$ by
$\Delta+\Pi_{\ker\Delta}$ and afterwards undo this finite-rank modification.
All assertions below concerning symbols are unchanged by this convention.

\subsection{From the model spectrum to compactness of a propagator}

\begin{lemma}\label{expanded:spectral-compactness}
Suppose the spectrum of a nonnegative self-adjoint operator $A$ with compact
resolvent is, outside a finite set, contained in a finite union of branches
\[
\lambda_{\alpha,k}=4(k+a_\alpha)(k+b_\alpha),\qquad
a_\alpha,b_\alpha\in\Z,
\]
where $k$ is a nonnegative integer and $\alpha$ belongs to a finite set.
Then $e^{2\pi i\sqrt A}-I$ is compact. The same conclusion holds for any
operator isospectral to $A$, with multiplicities included.
\end{lemma}

\begin{proof}
For fixed integers $a,b$, expansion at $k=\infty$ gives
\[
2\sqrt{(k+a)(k+b)}=2k+a+b+O(k^{-1}).
\]
The estimate is uniform over the finite collection of pairs $(a,b)$.
Consequently
\[
e^{2\pi i\sqrt{\lambda_{\alpha,k}}}-1=O(k^{-1}).
\]
If $\Pi_R$ is the spectral projection of $A$ onto eigenvalues at most $R$,
then $\Pi_R$ has finite rank and the spectral theorem gives
\[
\big\|(e^{2\pi i\sqrt A}-I)(I-\Pi_R)\big\|
=\sup_{\lambda\in\Spec(A),\ \lambda>R}
 |e^{2\pi i\sqrt\lambda}-1|\longrightarrow0.
\]
The operator is therefore a norm limit of finite-rank operators and is
compact. Increasing multiplicities cause no difficulty: the displayed
operator norm depends only on the largest scalar multiplier on each
eigenspace.

Two self-adjoint operators with the same discrete spectrum and the same
finite multiplicities admit a unitary intertwiner, obtained by choosing
orthonormal bases in corresponding eigenspaces. Functional calculus and
compactness are preserved by this unitary intertwiner. No geometric
identification of the underlying manifolds or vector bundles is being
assumed.
\end{proof}

For $(\CP^n,g_0)$ with holomorphic sectional curvature $4$, the required
finite collection of branches follows from Ikeda--Taniguchi,
Theorem~5.2(b) \cite{IkedaTaniguchi}, together with the Lefschetz decomposition into primitive
$(r,s)$-forms. In their normalization the branch polynomials have factors
with integral shifts; comparison with the scalar branch $k(k+n)$ changes
all eigenvalues by the common factor $4$ in our normalization. For each
fixed form degree only finitely many pairs of shifts occur. Thus the
previous lemma applies to every Hodge degree, not merely to functions.
Explicitly, before the common factor $4$, the polynomials in their table
are drawn from
\[
\begin{gathered}
(k+r)(k+n+2-s),\qquad (k+r)(k+n+1),\qquad k(k+n+1-s),\\
k(k+n),\qquad (k+n)(k+n+1),\qquad k(k+1),
\end{gathered}
\]
with the integral parameters and starting indices restricted as in that
theorem. The decomposition for any fixed form degree selects finitely
many such parameter values. Only the integrality of these shifts is
used here; no multiplicity formula is required for the compactness step.

\subsection{The microlocal input and the connection transport}

We state precisely the standard propagation facts used here. A
Laplace-type operator can be written uniquely as
\[
\Delta=\nabla^*\nabla+V,
\]
where $\nabla$ is a metric connection and $V$ is an endomorphism. Its
principal symbol is $|\xi|_g^2I$. The positive square root is an elliptic
pseudodifferential operator of order one with scalar principal symbol
$|\xi|_gI$. The half-wave propagator is a Fourier integral operator
associated with the graph of the homogeneous geodesic transformation.
We choose the sign convention for $G^t$ to agree with
$U(t)=e^{it\sqrt\Delta}$; changing this convention replaces $G^t$ by
$G^{-t}$ everywhere and does not change any periodicity conclusion.

The leading-symbol transport is parallel transport for $\nabla$. To see
why the potential $V$ does not create an additional leading phase, write
the symbol of $\sqrt\Delta$ as $a_1+a_0+a_{-1}+\cdots$. Comparing order
one terms in $(\sqrt\Delta)^2=\Delta$ determines $a_0$ from the
first-order connection part of $\Delta$; the zero-order term $V$ enters
only the equation for $a_{-1}$. In normal coordinates and a frame parallel
at a point, the covariant subprincipal symbol vanishes. In an arbitrary
frame its transport equation is precisely the parallel-transport equation
along the Hamiltonian trajectory. For the Hodge Laplacian the connection
is the induced Levi--Civita connection on exterior forms, by the
Weitzenb\"ock formula.

Working on half densities removes the coordinate divergence term from
scalar transport. The remaining scalar factor in the leading symbol is
the Maslov factor. At a time $L$ for which $G^L=I$, the canonical graph is
the diagonal, so $U(L)$ belongs to $\Psi^0$ and
\begin{equation}\label{expanded:return-symbol}
\sigma_0(U_p(L))(x,\xi)=\mu_L\,\Lambda^p H_L(x,\xi).
\end{equation}
Here $H_L$ denotes the induced return transport on covectors. Its tangent
version has the same eigenvalues and traces in the orthogonal situations
used below. The scalar phase $\mu_L$ is locally constant: it is the
Maslov index contribution for the continuously varying family of return
trajectories. Since $S^*X$ is connected when $X$ is connected and $D\geq2$,
it is a single constant. Its modulus is one. There is no additional
arbitrary scalar potential phase, by the subprincipal-symbol observation
above. These are the usual half-wave parametrix and transport statements;
see \cite{DG,HormanderIII,HormanderIV,Seeley}.
More precisely, Jakobson--Strohmaier \cite[Proposition 3.2]{JakobsonStrohmaier}
identify the subprincipal partial connection of $\sqrt{\nabla^*\nabla+V}$
with the lifted connection along the geodesic flow. This supplies the
bundle transport used in \eqref{expanded:return-symbol}.

\begin{lemma}\label{expanded:compact-return-implies-period}
For a Laplace-type operator, if $U(T)-I$ is compact on $L^2(X;E)$, then
$G^T=I$ on $S^*X$.
\end{lemma}

\begin{proof}
Put $K=U(T)-I$. If $A\in\Psi^0(X;E)$, then
\[
U(T)^*AU(T)-A=K^*A+AK+K^*AK
\]
is compact, because $A$ is bounded. Choose $A$ to have arbitrary scalar
principal symbol $a(x,\xi)I$. The bundle version of Egorov's theorem \cite[Theorem 1.4 and Proposition 3.3]{JakobsonStrohmaier} says
that the principal symbol of $U(T)^*AU(T)$ is $a\circ G^T$ times $I$:
the connection transports appearing on the two sides cancel because the
symbol is scalar. Thus the compact operator in the preceding display has
principal symbol $(a\circ G^T-a)I$.

An order-zero pseudodifferential operator that is compact on $L^2$ must
have zero principal symbol. One way to verify this last fact is to test
the operator on normalized high-frequency wave packets centered at a
point $(x_0,\xi_0)$ and a bundle vector $e$. These packets tend weakly to
zero, whereas the norms of their images tend to
$|\sigma_0(A)(x_0,\xi_0)e|$. Compactness makes those norms tend to zero,
so the principal symbol vanishes at every covector and bundle vector.
It follows that $a\circ G^T=a$ for every smooth function on $S^*X$.
Such functions separate points, hence $G^T=I$.
\end{proof}

Applying Lemmas~\ref{expanded:spectral-compactness} and
\ref{expanded:compact-return-implies-period} to the model spectrum gives
the common-period identity
\begin{equation}\label{period2}
G^{2\pi}=I.
\end{equation}

\subsection{The clean trace estimate with heat regularization}

\begin{lemma}\label{expanded:clean-heat-trace}
Let the fixed set of $G^L$ on $S^*X$ be clean, with connected components
$F_\alpha$ of dimensions $f_\alpha$. Then the contribution from
$F_\alpha$ to
\[
\Tr\bigl(U(L)e^{-\tau\Delta}\bigr),\qquad \tau\downarrow0,
\]
has order at most $\tau^{-(f_\alpha+1)/4}$. Contributions away from the
fixed set are bounded at $\tau=0$. In particular, if every
$f_\alpha<2D-1$, then
\[
\Tr\bigl(U(L)e^{-\tau\Delta}\bigr)=o(\tau^{-D/2}).
\]
If $G^L=I$ and $a_0$ is the principal symbol of $U(L)\in\Psi^0$, its
leading term is
\begin{equation}\label{expanded:identity-heat-leading}
(2\pi)^{-D}\int_{T^*X}
e^{-\tau|\xi|_g^2}\operatorname{tr}a_0(x,\xi)\,dx\,d\xi.
\end{equation}
In particular, if $\operatorname{tr}a_0$ is constant on $S^*X$, this is
\[
(4\pi\tau)^{-D/2}\Vol(X)\operatorname{tr}a_0.
\]
Changing $U(L)$ by an operator in $\Psi^{-1}$ changes its heat trace by
$O(\tau^{-(D-1)/2})$ for $D\geq2$, so it cannot change this leading
coefficient. This observation will later justify determining the constant
phase by comparison with the canonical K\"ahler involution.
\end{lemma}

\begin{proof}[Explanation of the clean stationary-phase calculation]
This is the fixed-time, heat-regularized form of the clean wave-trace
expansion in Duistermaat--Guillemin \cite[Theorem 4.5]{DG}.
For a finite-rank bundle with scalar principal symbol, the canonical
relation and stationary phase are unchanged: apply the scalar
oscillatory-integral calculation to the matrix-valued amplitude
entry by entry and then take the fibre trace. Connection transport
changes the leading amplitude, but not its order. We spell out the
Gaussian regularization below; no separate fixed-time exponent is
being assumed without calculation. The clean hypothesis means
\[
\ker(dG^L-I)_{\rho}=T_{\rho}\operatorname{Fix}(G^L)
\]
on the unit cosphere. Equivalently the canonical graph of $G^L$ and the
diagonal intersect cleanly in the nonzero cotangent bundle. Each fixed
component $F_\alpha$ becomes a cone with one additional radial variable.

After restricting a local oscillatory representation of the propagator
kernel to the diagonal, the stationary set is exactly this fixed cone.
The clean stationary-phase theorem integrates the nondegenerate
transverse variables. A component of cosphere dimension $f$ contributes
a radial integral whose leading amplitude has degree $(f-1)/2$:
\[
\int_1^\infty e^{-\tau\rho^2}
  \rho^{(f-1)/2}
  \bigl(c_0+O(\rho^{-1})\bigr)\,d\rho.
\]
The coefficient $c_0$ includes the trace of the bundle transport, the
Maslov factor and the clean stationary-phase density; it can vanish.
The degree $(f-1)/2$ is also the usual clean wave-trace order: taking the
Fourier transform in time produces a conormal singularity of order
$(f+1)/2$ at $L$. Introducing heat regularization replaces the radial
oscillatory cutoff by $e^{-\tau\rho^2}$ without changing the clean
critical set. Equivalently, set $h=\sqrt\tau$ and perform clean
semiclassical stationary phase, with principal heat amplitude
$e^{-|\xi|^2}$; this gives the same power of $h$.

The change of variables $u=\sqrt\tau\rho$ gives
\[
\int_0^\infty e^{-\tau\rho^2}\rho^{(f-1)/2}\,d\rho
=\frac12\Gamma\!\left(\frac{f+1}{4}\right)
 \tau^{-(f+1)/4}.
\]
Lower symbol terms have smaller growth; nonstationary pieces can be
integrated by parts in the high-frequency variables and are bounded.
This proves the stated order bound. Possible cancellation of leading
coefficients only improves the bound.

When the graph is the diagonal, no transverse stationary-phase reduction
is needed: the usual local pseudodifferential heat calculation directly
gives \eqref{expanded:identity-heat-leading}. Writing the symplectic
measure in metric cotangent coordinates and using
$\int_{\R^D}e^{-\tau|\xi|^2}\,d\xi=(\pi/\tau)^{D/2}$ proves the final
formula. For the identity return $f=2D-1$, the clean exponent is
$(f+1)/4=D/2$, in agreement with this direct calculation.
\end{proof}

For the intended application, cleanliness has a simple geometric proof,
so it is not an extra assumption about the unknown metric. If
$G^{2\pi}=I$, then $F=G^\pi$ is an involution on the compact manifold
$S^*X$. Average an arbitrary auxiliary Riemannian metric under $F$.
The map $F$ is then an isometry for this auxiliary metric. Around any fixed
point, the equivariance of its exponential map shows that the fixed set
is a smooth submanifold whose tangent space is the $+1$ eigenspace of
$dF$. This is precisely the clean condition. Its components are compact,
and only finitely many components occur. If a component had the full
dimension $2D-1$, then $F$ would agree with the identity on a nonempty
open set. An isometry of a connected Riemannian manifold is determined
by its value and derivative at one point, so $F$ would be the identity
everywhere. Consequently, unless $G^\pi=I$, all fixed components have
dimension less than $2D-1$ and the final little-oh estimate applies.

\section{Detailed topology and the reduction to total scalar curvature}
\label{sec:topology}

\begin{lemma}\label{expanded:compact-universal-cover}
If $G^T=I$ on the unit tangent bundle of a closed connected Riemannian
manifold of dimension at least two, its universal cover is compact.
\end{lemma}

\begin{proof}
Let $\widetilde X$ be the universal cover with the complete lifted
metric. The map $S\widetilde X\to SX$ is a connected regular covering,
with deck group $\pi_1(X)$. Since the lifted flow at time $T$ projects to
the identity, it is a deck transformation $\gamma$. Concretely, the
unique deck element carrying $\widetilde v$ to
$\widetilde G^T(\widetilde v)$ is locally constant in $\widetilde v$;
connectedness makes it one element for every $\widetilde v$.

Let $R(v)=-v$ denote reversal. It commutes with every lifted deck action,
and the reversibility of geodesic flow gives
$R\widetilde G^TR=\widetilde G^{-T}$. Thus $\gamma=\gamma^{-1}$ and
$\widetilde G^{2T}=I$. Every unit-speed lifted geodesic has period
$2T$, possibly a smaller prime period. A minimizing geodesic segment
has length at most $T$: if its length is between $T$ and $2T$, the
oppositely oriented part of the same closed geodesic is shorter; a
segment of length at least $2T$ contains a nontrivial loop and cannot be
minimizing. Hopf--Rinow supplies minimizing geodesics between all pairs
of points, so $\operatorname{diam}\widetilde X\leq T$.
The complete bounded manifold $\widetilde X$ is compact by Hopf--Rinow.
\end{proof}

\begin{proposition}\label{expanded:projective-topology}
A compact connected K\"ahler manifold whose geodesic flow has a common
period is simply connected and has integral cohomology ring
$\Z[x]/(x^{n+1})$, with $|x|=2$.
\end{proposition}

\begin{proof}
The previous lemma makes its simply connected universal cover compact.
The Bott--Samelson theorem together with McCleary's refinement says that
a compact simply connected Besse manifold has the integral cohomology
ring of a compact rank-one symmetric space. We use precisely this
integral conclusion. It is stated explicitly in the introduction of
Radeschi--Wilking \cite{RadeschiWilking}, which cites
Besse, Theorem~7.37 \cite{Besse}, and McCleary, Corollary~A \cite{McCleary}.

The K\"ahler form lifts to a K\"ahler form on the cover, and its degree-two
class is nonzero. Indeed, if it were exact, Stokes' theorem would make
its positive volume integral $\int\widetilde\omega^n$ vanish. The sphere
in real dimension above two, the quaternionic projective spaces, and the
Cayley plane have zero real degree-two cohomology. The only remaining
rank-one model is $\CP^n$, including $S^2=\CP^1$. Hence the universal
cover has the asserted integral ring.

Every deck transformation preserves $\widetilde\omega$ exactly, since it
is pulled back from the base. Therefore it fixes the generator of
$H^2(\widetilde M;\Q)$ and its powers, and acts trivially on all
rational cohomology. Its Lefschetz number is
\[
\sum_j(-1)^j\Tr\bigl(\gamma^*:H^j(\widetilde M;\Q)\bigr)=n+1.
\]
The Lefschetz fixed-point theorem \cite{Hatcher} implies that it has a fixed point.
Nonidentity deck transformations of a covering have none. The deck
group is consequently trivial, and therefore
\begin{equation}\label{topology}
\pi_1(M)=0,\qquad H^*(M;\Z)=\Z[x]/(x^{n+1}),\quad |x|=2.
\end{equation}
\end{proof}

This proposition does not assert that $M$ is homeomorphic to projective
space, and no such implication is needed. In particular, the subsequent
argument does not identify Pontryagin classes from a homotopy type.

Choose the sign of the integral generator so that $[\omega]=t x$ with
$t>0$. The integral ring and positivity imply $\int_Mx^n=1$.
Consequently
\[
\Vol(M,g)=\frac1{n!}\int_M\omega^n=\frac{t^n}{n!}.
\]
For the normalized model this volume is $\pi^n/n!$. Equality of the
leading heat coefficients therefore gives $t=\pi$, without rescaling
the unknown metric after its spectrum has been fixed.

Write $c_1(M)=kx$. Let $\rho$ be the ordinary real Ricci form, so that
$[\rho]=2\pi c_1(M)$ and $s=2\Lambda_\omega\rho$. The pointwise identity
$\rho\wedge\omega^{n-1}=(s/2n)\omega^n$ gives
\[
\int_Ms\,dV
=\frac{4\pi}{(n-1)!}\,
c_1(M)[\omega]^{n-1}
=\frac{4\pi k\pi^{n-1}}{(n-1)!}.
\]
Dividing by $\pi^n/n!$ yields
\begin{equation}\label{meanscalar}
\overline s:=\frac1V\int_Ms\,dV=4nk.
\end{equation}
Thus recovering the model's total scalar curvature is exactly the same
as establishing $k=n+1$.

Finally, $x=[\omega]/\pi$ is an integral class of type $(1,1)$. The
Lefschetz $(1,1)$ theorem \cite{Huybrechts} supplies a holomorphic line bundle $L$ with
$c_1(L)=x$. Its class contains the positive form $\omega/\pi$; by the
$\partial\bar\partial$ lemma one can choose a Hermitian metric on $L$
with positive curvature. Kodaira's theorem \cite{Kodaira} then makes $L$ ample and $M$
projective. If $k=n+1$, the anticanonical class is positive, so $M$ is
Fano. Moreover $H^1(M;\C)=0$, hence $H^1(M,\mathcal O)=0$; the
exponential sequence makes $c_1:\operatorname{Pic}(M)\to H^2(M;\Z)$
injective. Therefore $-K_M\cong L^{n+1}$, not merely numerically.
The Kobayashi--Ochiai theorem \cite{KobayashiOchiai} gives $M\cong\CP^n$ biholomorphically.
Under this identification $x$ is the positive hyperplane generator and
$[\omega]=\pi x$ is precisely the normalized Fubini--Study class.

\subsection{The ordinary first coefficient and the exceptional equation}

The integrated Hodge heat formula \cite{Patodi,PatodiEigenforms,Gilkey} is
\[
 a_{1,p}=\beta_p\int_Ms\,dV,\qquad
 \beta_p=\frac{N}{6}-\binom{D-2}{p-1}.
\]
The binomial ratio, with the zero convention when $p=0$, is
\[
 \frac{\binom{D-2}{p-1}}{\binom Dp}
   =\frac{p(D-p)}{D(D-1)}.
\]
It follows from $D=2n$ that
\[
 \beta_p=\frac{N}{D(D-1)}
       \left(\frac{D(D-1)}6-p(D-p)\right)
       =\frac{Nq}{D(D-1)}.
\]
If $q\ne0$, equality of $a_{1,p}$ recovers total scalar curvature;
the preceding argument identifies the complex manifold and its
K\"ahler class. Such pairs can proceed directly to Section~\ref{sec:metric}.

The remaining sections address precisely
\begin{equation}\label{exception}
 q=0\quad\Longleftrightarrow\quad
 p(D-p)=\frac{D(D-1)}6.
\end{equation}
For $n\geq1$ this equation excludes $p=0$ and $p=n$.
In dimensions $n=1,2$ there is no remaining integer solution;
therefore every exceptional pair has $0<p<n$ and $n\geq3$.
The Pell calculation below lists all such pairs; the first three are
$(p,n)=(1,3),(20,48),(285,675)$.

\section{The exceptional spectrum determines return holonomy}
\label{sec:return}

\subsection{The model's return holonomy and why its leading trace survives}

The holonomy assertion for Fubini--Study space can be seen directly from
the Hopf model. Let $e_0,e_1,\ldots,e_n$ be the standard unitary basis of
$\C^{n+1}$, and consider the horizontal great circle
\[
u(t)=\cos t\,e_0+\sin t\,e_1\in S^{2n+1},\qquad 0\leq t\leq\pi.
\]
Its image under the Hopf projection is a unit-speed geodesic for the
Fubini--Study metric of holomorphic sectional curvature $4$. The
horizontal spaces are the complex orthogonal complements of $u(t)$.
The fields $\dot u(t)$ and $i\dot u(t)$ project to parallel fields along
the geodesic: the derivative of the first is normal to the sphere, and
the horizontal part of the derivative of the second vanishes. The
constant fields $e_j,ie_j$, $j\geq2$, are horizontal and project to
parallel fields as well.

At $t=\pi$ the lift is $-e_0$. Identifying this endpoint with the initial
Hopf representative $e_0$ multiplies horizontal tangent vectors by
$-1$. Since $\dot u(\pi)=-e_1$, this endpoint identification preserves
the projected geodesic tangent and its $J$-multiple, while it negates
every projected field $e_j,ie_j$ with $j\geq2$. Therefore
\[
H_\pi=I\quad\hbox{on }\operatorname{span}_{\R}(v,Jv),\qquad
H_\pi=-I\quad\hbox{normally}.
\]
The exterior character is
\[
\operatorname{tr}\Lambda^pH_\pi
=[z^p](1+z)^2(1-z)^{D-2}
=(-1)^p\left(\binom Dp-4\binom{D-2}{p-1}\right).
\]
At an exceptional pair, the binomial coefficient on the right is
$N/6$, so the character is $(-1)^pN/3\ne0$. The scalar Maslov factor
has modulus one and is constant; hence it cannot cancel this coefficient
after integration over $S^*\CP^n$. Formula
\eqref{expanded:identity-heat-leading} gives a nonzero leading
$\tau^{-D/2}$ term. Ordinary isospectrality makes the entire
heat-regularized half-wave trace equal on the two manifolds, so the
preceding clean-trace estimate forces the shorter common return time
\begin{equation}\label{periodpi}
G^\pi=I
\end{equation}
on the unknown manifold.

\subsection{The involutivity of the unknown return holonomy}

Once $G^\pi=I$, write the symbol as $\mu\Lambda^pH(v)$ using
\eqref{expanded:return-symbol}. The symbol calculus and
$U(2\pi)=U(\pi)^2$ give
\[
\sigma_0(U(2\pi))=\mu^2\Lambda^p(H(v)^2).
\]
Its difference from $I$ is compact, so this symbol is $I$.
For completeness, if an orthogonal map $A$ satisfies
$\Lambda^pA=cI$ with $0<p<D$, diagonalize $A$ over $\C$ with
eigenvalues $\alpha_1,\ldots,\alpha_D$. Given $i\ne j$, choose a set
of $p-1$ indices disjoint from $i,j$; this is possible because
$p\leq D-1$. Comparing the two $p$-fold eigenvalue products gives
$\alpha_i=\alpha_j$. Thus $A=\alpha I$. Apply this to $A=H(v)^2$.
Since parallel transport around the returned geodesic fixes its tangent,
$A(v)=v$, whence $\alpha=1$. We obtain
\begin{equation}\label{involhol}
H(v)^2=I,\qquad \mu^2=1.
\end{equation}
As $\nabla J=0$, each $H(v)$ is complex linear. Its negative complex
eigenspace rank is locally constant and hence constant on $SM$.
Denote the negative complex eigenspace rank by $r$.
Since $v,Jv$ are fixed, $0\leq r\leq n-1$. On the real tangent space
there are $2r$ negative and $D-2r$ positive eigenvalues, so
\[
\operatorname{tr}\Lambda^pH(v)
       =[z^p](1-z)^{2r}(1+z)^{D-2r}.
\]
This character is independent of $v$. The identity-return heat formula,
equal volumes, and equality of the ordinary regularized wave traces give
\[
 \mu_M[z^p](1-z)^{2r}(1+z)^{D-2r}
       =\mu_0\frac{(-1)^pN}{3}.
\]
Taking absolute values eliminates the unknown unit phases and yields
\begin{equation}\label{character}
 \left|[z^p](1-z)^{2r}(1+z)^{D-2r}\right|=\frac{N}{3}.
\end{equation}
We next show that this forces the full model holonomy.

\section{A uniform Krawtchouk estimate}

The generating function, duality, recurrence, and orthogonality relations
for Krawtchouk polynomials are classical; see
\cite[\S9.11]{KoekoekLeskySwarttouw} and
\cite[\S\S18.19, 18.22--18.23]{DLMF}.
We give a self-contained proof of the particular uniform estimate required
in \eqref{character}, including the normalization and the elementary
estimates. Throughout this section,
\[
 D=2n,\qquad 0<p<n,\qquad
 p(D-p)=\frac{D(D-1)}6,
 \qquad S=D-2p>0.
\]
The exceptional equation is equivalently
\begin{equation}\label{expanded-S-square}
 S^2=\frac{D(D+2)}3.
\end{equation}
We use the convention that a binomial coefficient is zero when its lower
index lies outside the interval from zero to its upper index.

\begin{lemma}\label{kraw}
Suppose $D=2n$ and $p$ satisfy the preceding conditions.  If
$0\le r\le n-1$ and
\[
 \left|[t^p](1-t)^{2r}(1+t)^{D-2r}\right|
 =\frac13\binom Dp,
\]
then $r=1$ or $r=n-1$.
\end{lemma}

\begin{proof}
\emph{Normalization and duality.}
Define the binary Krawtchouk coefficients by
\[
 K_j(a;D)=[z^j](1-z)^a(1+z)^{D-a}
 =\sum_{\ell}(-1)^\ell\binom a\ell\binom{D-a}{j-\ell}.
\]
For every admissible summation index,
\[
 \frac{\binom p\ell\binom{D-p}{j-\ell}}{\binom Dj}
 =
 \frac{\binom j\ell\binom{D-j}{p-\ell}}{\binom Dp}.
\]
Both sides count the same intersection-size probability for a $p$-subset
and a $j$-subset of a $D$-element set; the equality also follows immediately
by expanding the four factorials.  Summing with the sign $(-1)^\ell$ gives
\begin{equation}\label{expanded-kraw-duality}
 f_j:=\frac{K_j(p;D)}{\binom Dj}
 =\frac{K_p(j;D)}{\binom Dp}.
\end{equation}
The desired conclusion is therefore a statement about the even indices
$j=2r$ for which $|f_j|=1/3$.

\emph{The recurrence.}
Let $F(z)=(1-z)^p(1+z)^{D-p}$.  Its logarithmic derivative gives
\[
 (1-z^2)F'(z)=(S-Dz)F(z).
\]
Comparing coefficients of $z^j$ gives
\[
 (j+1)K_{j+1}(p;D)
 =S K_j(p;D)-(D-j+1)K_{j-1}(p;D).
\]
Using
$\binom D{j+1}/\binom Dj=(D-j)/(j+1)$ and
$\binom D{j-1}/\binom Dj=j/(D-j+1)$, we obtain
\begin{equation}\label{expanded-kraw-recurrence}
 (D-j)f_{j+1}=S f_j-jf_{j-1}.
\end{equation}
Here $f_0=1$ and $f_1=S/D$.  Substitution of
\eqref{expanded-S-square} yields
\[
 f_2=\frac{S^2-D}{D(D-1)}=\frac13,
 \qquad
 f_3=\frac{S(D-6)}{3D(D-2)}.
\]
For $D>6$, we have $0<S<D$ and $0<D-6<D-2$, hence
$0<f_3<1/3$.

For $3\le j\le p$, the identity $S=D-2p$ implies
\[
 \frac{S+j}{D-j}\le1.
\]
Starting with $|f_2|=1/3$ and $|f_3|<1/3$, the triangle inequality in
\eqref{expanded-kraw-recurrence} gives successively
\begin{equation}\label{smallk}
 |f_j|<\frac13\qquad(3\le j\le p+1).
\end{equation}
The first induction step is strict because the coefficient $S$ of $f_3$
is positive; every subsequent step involves two already strict bounds.

\emph{Orthogonality and the middle range.}
The generating-function identity
\begin{align*}
 &\sum_{j=0}^D\binom Dj
 (1-z)^j(1+z)^{D-j}(1-w)^j(1+w)^{D-j}\\
 &\hspace{15mm}=
 \big((1-z)(1-w)+(1+z)(1+w)\big)^D
 =2^D(1+zw)^D
\end{align*}
implies, on comparing the coefficient of $z^p w^q$,
\[
 \sum_{j=0}^D\binom Dj K_p(j;D)K_q(j;D)
 =2^D\binom Dp\,\delta_{pq}.
\]
Thus \eqref{expanded-kraw-duality} gives
\begin{equation}\label{expanded-kraw-orthogonality}
 \sum_{j=0}^D\binom Dj f_j^2
 =\frac{2^D}{\binom Dp}.
\end{equation}
Each summand is nonnegative.  For $p\le j\le D-p$, unimodality and
symmetry of the binomial coefficients give $\binom Dj\ge\binom Dp$,
so
\begin{equation}\label{middlek}
 |f_j|\le\sqrt{\frac{2^D}{\binom Dp\binom Dj}}
 \le\frac{2^{D/2}}{\binom Dp}.
\end{equation}

\emph{The first exceptional dimensions.}
We justify the dimension gap used to estimate the last expression.
By \eqref{expanded-S-square}, the positive integers $X=D+1$ and $Y=S$
satisfy $X^2-3Y^2=1$, with $X$ odd and $Y$ even.  The positive solutions
of this Pell equation are
\[
 X+Y\sqrt3=(2+\sqrt3)^b,\qquad b\ge1.
\]
Here is the elementary descent proof of the assertion being used.  If
$Y\ge2$, multiplication by $2-\sqrt3$ replaces $(X,Y)$ by
\[
 (X',Y')=(2X-3Y,\,2Y-X).
\]
The equality $X^2=3Y^2+1$ gives
$\sqrt3Y<X<2Y$ for $Y>1$.  Therefore $X'>0$ and $0<Y'<Y$.
The new pair again satisfies $X'^2-3Y'^2=1$.  Repeated descent reaches
$(X,Y)=(2,1)$, which maps to $(1,0)$.  Reversing the descent proves the
displayed parametrization.  Modulo two, the recurrence for multiplication
by $2+\sqrt3$ interchanges the two coordinate parities.  Thus $X$ odd and
$Y$ even requires $b$ even, and
\[
 D+1+S\sqrt3=(7+4\sqrt3)^a,\qquad a\ge1.
\]
The first two terms are $7+4\sqrt3$ and $97+56\sqrt3$.
They give $(D,p)=(6,1)$ and $(96,20)$.  The coordinates strictly increase
with $a$, so every other exceptional pair has $D\ge96$.

\emph{An explicit binomial lower bound.}
Suppose now that $D\ge96$ and put $\alpha=p/D$.  Then
\[
 \alpha=\frac12\left(1-\sqrt{\frac{D+2}{3D}}\right),
 \qquad \frac15<\alpha<\frac12.
\]
The first strict inequality follows already from $D>25$, since
$(D+2)/(3D)<9/25$ is equivalent to $D>25$.
Consider the binomial probability distribution with $D$ trials and
success probability $\alpha=p/D$.  The atom at $p$ is a mode: the ratio
of the atom at $j+1$ to that at $j$ is
\[
 \frac{D-j}{j+1}\frac{p}{D-p},
\]
which is greater than one at $j=p-1$ and less than one at $j=p$.
The largest of $D+1$ nonnegative numbers summing to one is at least
$1/(D+1)$.  Consequently
\[
 \binom Dp\alpha^p(1-\alpha)^{D-p}\ge\frac1{D+1},
 \qquad
 \binom Dp\ge\frac{2^{D\mathcal H(\alpha)}}{D+1},
\]
where
$\mathcal H(\alpha)=-\alpha\log_2\alpha-(1-\alpha)\log_2(1-\alpha)$.
Since $\mathcal H'(\alpha)=\log_2((1-\alpha)/\alpha)>0$ on $(0,1/2)$,
\[
 \mathcal H(\alpha)>\mathcal H(1/5)
 =\log_2 5-\frac85>\frac7{10}.
\]
The last inequality is the integer comparison
$5^{10}=9765625>8388608=2^{23}$.
Also $2^{D/5}>3(D+1)$ for $D\ge96$: the quotient
$2^{D/5}/(D+1)$ has positive logarithmic derivative
$(\log 2)/5-1/(D+1)$ on this interval, and at $D=96$ we have
$2^{96/5}>2^{19}>291=3\cdot97$.
It follows that
\[
 \binom Dp>
 \frac{2^{7D/10}}{D+1}
 >3\,2^{D/2}.
\]
Equation \eqref{middlek} is therefore strictly less than $1/3$ throughout
the middle range.

\emph{The large indices and the conclusion.}
The generating function gives
\[
 z^D F(1/z)=(-1)^pF(z),
 \qquad f_{D-j}=(-1)^p f_j.
\]
Combining this symmetry, \eqref{smallk}, and the middle-range estimate
proves
\[
 |f_j|<1/3\qquad(3\le j\le D-3).
\]
Among the even indices $j=2r$ with $0\le r\le n-1$, the remaining
possibilities are $0,2,D-2$.  Since $f_0=1$, $f_2=1/3$, and
$f_{D-2}=(-1)^p/3$, exactly $j=2,D-2$ are possible.  In the single small
case $(D,p)=(6,1)$, direct calculation gives
$K_1(2r;6)=6-4r$ for $r=0,1,2$; the equation $|6-4r|=2$ again gives
$r=1,2=n-1$.  This proves the lemma.
\end{proof}

\section{A mod-two obstruction eliminates the complementary holonomy}

We spell out the bundle construction, the cohomology classes, and the
two coefficients that yield the obstruction.  Write $c_1(TM,J)=kx$,
where $x$ is the integral generator chosen in \eqref{topology}; the same
letter $x$ will denote its reduction in mod-two cohomology and its
pullback to any bundle over $M$.

\begin{lemma}\label{lineobstruction}
Assume \eqref{topology}, \eqref{exception}, \eqref{periodpi}, and
\eqref{involhol}.  The negative complex eigenspace of $H_\pi(v)$ cannot
have complex dimension one.
\end{lemma}

\begin{proof}
\emph{Smoothness and descent of the eigenline.}
Suppose for contradiction that this dimension is one.  The geodesic
flow and the ordinary differential equation for parallel transport
depend smoothly on the initial unit vector $v$.  Since $H(v)^2=I$, the
operator $(I-H(v))/2$ is a smooth family of projections of constant
complex rank one.  Its image is therefore a smooth complex line bundle
over $SM$.

The geodesic with initial velocity $-v$ traverses the closed geodesic
with initial velocity $v$ in the opposite direction.  Parallel transport
along a reversed closed curve is the inverse transport.  Accordingly
\[
 H(-v)=H(v)^{-1}=H(v).
\]
The eigenlines over $v$ and $-v$ are thus the very same subspace of
$T_{\pi(v)}M$, with the identity identification.  Quotienting by this
free antipodal action produces a well-defined smooth complex line
bundle $E$ over
\[
 P=P_\R(TM),\qquad \pi:P\longrightarrow M.
\]
This descent uses only $v\mapsto-v$; no invariance under
$v\mapsto e^{i\theta}v$ is assumed.

\emph{A second, tautological complex line.}
Let $\tau\subset\pi^*(TM)_\R$ be the real tautological line bundle.
Its complexification $\eta=\tau\otimes_\R\C$ has a natural complex-linear
embedding into $\pi^*(TM,J)$ given by
\[
 v\otimes(a+ib)\longmapsto av+bJv.
\]
This is well-defined by real bilinearity, is complex-linear because
multiplication by $i$ in the target is $J$, and is injective because
$v$ and $Jv$ are linearly independent over $\R$ for $v\ne0$.
Parallel transport fixes the velocity of a geodesic, and $J$ is parallel
on a K\"ahler manifold.  Thus $H(v)$ fixes both $v$ and $Jv$.
An orthogonal involution has orthogonal positive and negative eigenspaces,
so $\eta$ and $E$ are orthogonal subbundles of $\pi^*TM$.
Consequently
\[
 F=\pi^*TM/(\eta\oplus E)
\]
is a complex vector bundle of rank $n-2$.

\emph{The necessary cohomology of $P$.}
All classes for the rest of this proof have coefficients in $\F_2$.
Set
\[
 z=w_1(\tau)\in H^1(P;\F_2),\qquad u=z^2\in H^2(P;\F_2).
\]
On each fibre $\R P^{2n-1}$, $z$ restricts to the standard degree-one
generator.  Since
$H^*(\R P^{2n-1};\F_2)=\F_2[z]/(z^{2n})$, Leray--Hirsch gives the
free $H^*(M;\F_2)$-module decomposition
\begin{equation}\label{expanded-projective-module}
 H^*(P;\F_2)
 =\bigoplus_{j=0}^{2n-1}\pi^*H^*(M;\F_2)\,z^j.
\end{equation}
This is the part of the real projective-bundle theorem that is needed;
the relation involving $z^{2n}$ and the Stiefel--Whitney classes of $TM$
does not enter the proof.  In particular, in real degree $2n-2$ the
monomials
\[
 x^a u^{n-1-a}=x^a z^{2n-2-2a},\qquad 0\le a\le n-1,
\]
are linearly independent over $\F_2$.

The underlying real bundle of $\eta$ is $\tau\oplus\tau$.  Therefore
\[
 w(\eta_\R)=(1+z)^2=1+z^2,
 \qquad c_1(\eta)\bmod2=w_2(\eta_\R)=u.
\]
Here we used the standard identity $c_1(L)\bmod2=w_2(L_\R)$ for a
complex line bundle.  Because $H^1(M;\F_2)=0$, the degree-two part of
\eqref{expanded-projective-module} has basis $x,u$.  Hence for uniquely
determined $a,b\in\F_2$,
\[
 c_1(E)\bmod2=ax+bu.
\]

\emph{The two independent Chern-class coefficients.}
Whitney multiplicativity in the direct sum
$\pi^*TM\cong\eta\oplus E\oplus F$ gives
\[
 c(F)=\frac{\pi^*c(TM)}{(1+u)(1+ax+bu)}.
\]
Treat $x$ and $u$ as having complex degree one, or real degree two.
Since $F$ has rank $n-2$, its component of complex degree $n-1$ is zero:
\begin{equation}\label{expanded-Chern-zero}
 0=c_{n-1}(F)
 =\left[\frac{\pi^*c(TM)}{(1+u)(1+ax+bu)}\right]_{n-1}.
\end{equation}

The coefficient of $u^{n-1}$ uses only the constant term of $c(TM)$
and the terms with no $x$ in the denominators.  It is
\[
 [u^{n-1}]\frac1{(1+u)(1+bu)}
 =\sum_{j=0}^{n-1}b^j.
\]
For $b=0$ this equals one.  For $b=1$ it equals $n$ modulo two.
Independence of the monomials above shows that \eqref{expanded-Chern-zero}
is impossible when $n$ is odd.  If $n$ is even, it forces $b=1$.

Assume therefore $n$ even and $b=1$.  Up to terms divisible by $x^2$, in
characteristic two we have
\begin{align*}
 \frac{\pi^*c(TM)}{(1+u)(1+ax+u)}
 &=\frac{1+kx+O(x^2)}{(1+u)^2}
   +\frac{ax}{(1+u)^3}+O(x^2)\\
 &=\frac1{(1+u)^2}
   +x\left(\frac{k}{(1+u)^2}+\frac{a}{(1+u)^3}\right)+O(x^2).
\end{align*}
The formal power-series identity
\[
 (1+u)^{-r}
 =\sum_{j\ge0}\binom{r+j-1}{j}u^j
 \qquad\text{over }\F_2
\]
follows by multiplying $r$ copies of the geometric series and counting
compositions of $j$ into $r$ nonnegative integers.  Thus the coefficient
of $xu^{n-2}$ in \eqref{expanded-Chern-zero} is
\begin{equation}\label{expanded-second-coefficient}
 k(n-1)+a\binom n2\pmod2.
\end{equation}

\emph{Parity forced by the exceptional equation and Wu's identity.}
Rewriting \eqref{exception} in the original variables gives
\[
 n(n+1)=3(n-p)^2.
\]
Modulo four the left side is zero for $n\equiv0,3\pmod4$ and two for
$n\equiv1,2\pmod4$, whereas the right side is zero or three.  Therefore
$n\equiv0$ or $3\pmod4$.  In the even case currently under consideration,
$n\equiv0\pmod4$, so $n-1$ is odd and
$\binom n2=(n/2)(n-1)$ is even.  Equation
\eqref{expanded-second-coefficient} consequently forces $k$ even.

We now compute the actual parity of $k$ without making any assertion
about Pontryagin classes.  The complex orientation gives $w_1(TM)=0$.
Wu's identity says that the degree-two Wu class $v_2$ satisfies
\[
 \langle v_2\smile\alpha,[M]\rangle
 =\langle\operatorname{Sq}^2\alpha,[M]\rangle,
 \qquad \alpha\in H^{2n-2}(M;\F_2),
\]
and the relation $w=\operatorname{Sq}(v)$ gives $v_2=w_2$ on an oriented
manifold.  Take $\alpha=x^{n-1}$.  Because $x$ lifts to an integral
class, $\operatorname{Sq}^1x=0$; and by the unstable square identity,
$\operatorname{Sq}^2x=x^2$.  The Cartan formula therefore gives
\[
 \operatorname{Sq}^2(x^{n-1})=(n-1)x^n.
\]
Since $\langle x^n,[M]\rangle=1$ and $H^2(M;\F_2)$ is generated by $x$,
\[
 w_2(TM)=(n-1)x=(n+1)x.
\]
For a complex vector bundle, reduction of $c_1$ modulo two is $w_2$ of
the underlying real bundle (by the line-bundle identity and the
splitting principle).  Hence
\[
 k\equiv n+1\pmod2.
\]
For the even $n$ under consideration, $k$ is odd, contradicting the
even parity forced by \eqref{expanded-second-coefficient}.
Both parities of $n$ have led to contradictions, so the negative
eigenspace cannot be a complex line.
\end{proof}

Combining Lemmas~\ref{kraw} and \ref{lineobstruction}, the negative
complex eigenspace has dimension $n-1$.  Its orthogonal positive
eigenspace has complex dimension one and contains the complex line
spanned by $v$.  Thus the parallel transport is determined, not merely
its trace:
\begin{equation}\label{modelhol}
 H(v)=+I\text{ on }\operatorname{span}_\R(v,Jv),\qquad
 H(v)=-I\text{ on its orthogonal complement}.
\end{equation}

The projective-bundle theorem, Whitney multiplicativity, and the Wu
identities used above are the standard characteristic-class results
in \cite{MilnorStasheff}. The Krawtchouk identities and estimates have
been derived directly from the generating function.

\section{An exact K\"ahler involution and its weighted heat coefficient}

\subsection{K\"ahler anticommutation identities and the Green operator}

All operators in this section act on the complexified exterior algebra
$\Omega^*(M;\C)$, with its $L^2$ inner product.  Put
\[
 d=\partial+\bar\partial,
 \qquad d^c=i(\bar\partial-\partial),
 \qquad \delta=d^*,\qquad \delta^c=(d^c)^*.
\]
We write $\{A,B\}=AB+BA$ for the anticommutator of two odd operators.
The K\"ahler identities imply
\[
 \{\partial,\partial^*\}
 =\{\bar\partial,\bar\partial^*\}
 =\tfrac12\Delta,
 \qquad
 \{\partial,\bar\partial^*\}
 =\{\bar\partial,\partial^*\}=0.
\]
Together with $\partial^2=\bar\partial^2=0$ and
$\partial\bar\partial+\bar\partial\partial=0$, these give
\begin{align}
 &d^2=(d^c)^2=\delta^2=(\delta^c)^2=0,\label{carzero}\\
 &\{d,\delta\}=\{d^c,\delta^c\}=\Delta,\label{carnonzero}\\
 &\{d,d^c\}=\{d,\delta^c\}
  =\{\delta,d^c\}=\{\delta,\delta^c\}=0.\label{carcross}
\end{align}
For example,
$\delta^c=i(\partial^*-\bar\partial^*)$, so
$\{d,\delta^c\}=i(\Delta_\partial-\Delta_{\bar\partial})=0$.
The other mixed identities follow in the same way or by taking adjoints.
Each of these four odd operators commutes with $\Delta$ and preserves every
positive eigenspace of the total Hodge Laplacian, while changing form degree
by one.

Let $\Pi$ be orthogonal projection onto the harmonic forms, and define the
Green operator $G$ by
\[
 G|_{\ker\Delta}=0,
 \qquad G|_{E_\lambda}=\lambda^{-1}I\quad(\lambda>0).
\]
Thus $\Delta G=G\Delta=I-\Pi$.  Elliptic regularity makes $\Pi$ a smoothing
operator, and elliptic pseudodifferential calculus gives
$G\in\Psi^{-2}$, with principal symbol $|\xi|^{-2}I$.
Here is the usual parametrix justification.  If $Q\in\Psi^{-2}$ is a
parametrix with $\Delta Q=I-S$, $S\in\Psi^{-\infty}$, then
\[
 G-Q=GS-\Pi Q\in\Psi^{-\infty}.
\]
The first term is smoothing because $G$ improves Sobolev order by two at
every order, and the second has smooth finite-dimensional range.  This also
shows that adding or removing the harmonic projection never changes a
principal symbol.

The K\"ahler identities used here are the standard identities in
Huybrechts \cite{Huybrechts}.
The parametrix argument is the ordinary elliptic calculus of
H\"ormander \cite[Chapter XVIII]{HormanderIII}.  The displayed calculations specify the conventions needed
here, so no sign convention for $d^c$ from either reference is implicit.

Define degree-preserving operators
\[
 P=d\delta G,\qquad P^c=d^c\delta^cG,
 \qquad \mathcal K=(I-2P)(I-2P^c).
\]
They are bounded order-zero pseudodifferential operators.
On a positive eigenspace introduce
\[
 a=\lambda^{-1/2}d,\qquad b=\lambda^{-1/2}d^c,
 \qquad a^*=\lambda^{-1/2}\delta,
 \qquad b^*=\lambda^{-1/2}\delta^c.
\]
Equations \eqref{carzero}--\eqref{carcross} become the canonical
anticommutation relations
\[
 a^2=b^2=0,\quad aa^*+a^*a=bb^*+b^*b=I,
 \quad ab=-ba,\quad ab^*=-b^*a,
\]
and their adjoints.  In particular $P=aa^*$ and $P^c=bb^*$ are
self-adjoint projections.  For example,
\[
 (aa^*)^2=a(I-aa^*)a^*=aa^*.
\]
They commute: moving $a,a^*$ past $b,b^*$ introduces two minus signs.
On harmonic forms both projections vanish.  Therefore, globally,
\begin{equation}\label{Kexact}
 \mathcal K^*=\mathcal K,\qquad \mathcal K^2=I,
 \qquad [\mathcal K,\Delta]=0.
\end{equation}

\subsection{Principal symbol and its relation to the wave operator}

For a real covector $\xi\ne0$, write $e(\xi)$ for exterior multiplication
by $\xi$ and $\iota_{\xi^\sharp}$ for contraction.  The symbols of $d$ and
$\delta$ are $i e(\xi)$ and $-i\iota_{\xi^\sharp}$, respectively.  Hence
\[
 \sigma_0(P)(x,\xi)
 =|\xi|^{-2}e(\xi)\iota_{\xi^\sharp}.
\]
This is the orthogonal projection onto the exterior forms that contain the
factor $\xi/|\xi|$.  Similarly $\sigma_0(P^c)$ is the projection onto forms
containing the orthogonal factor $J\xi/|\xi|$; the possible sign in the
symbol of $d^c$ cancels against its adjoint.  Choose an orthonormal coframe
whose first two elements are $\xi/|\xi|,J\xi/|\xi|$.  On each exterior
monomial, $\sigma_0(\mathcal K)$ is multiplication by $(-1)^\ell$, where
$\ell$ is the number of these two factors occurring in the monomial.
Thus $\sigma_0(\mathcal K)$ is precisely the exterior action of the
reflection that is negative on their two-plane and positive on its
orthogonal complement.

The holonomy in \eqref{modelhol} is the opposite reflection.  On degree
$p$, changing a linear map to its negative multiplies its exterior action
by $(-1)^p$.  Consequently
\[
 \Lambda^pH=(-1)^p\sigma_0(\mathcal K_p).
\]
The principal-symbol formula for the half-wave operator therefore implies
\begin{equation}\label{waveK}
 U_p(\pi)-\varepsilon\mathcal K_p\in\Psi^{-1}
\end{equation}
for a constant phase $\varepsilon$ of absolute value one.  This phase
includes both the Maslov factor and the displayed factor $(-1)^p$.

\subsection{The four summands in each positive eigenspace}

Fix $\lambda>0$.  Let $E_\lambda^j$ be its degree-$j$ eigenspace and put
\[
 B_\lambda^j=E_\lambda^j\cap\ker a^*\cap\ker b^*,
 \qquad b_j(\lambda)=\dim_\C B_\lambda^j,
 \qquad m_j(\lambda)=\dim_\C E_\lambda^j.
\]
Degrees outside $0,\ldots,D$ are interpreted as zero.  The two commuting
projections $P,P^c$ split $E_\lambda^j$ orthogonally into their four joint
eigenspaces, with projection eigenvalues $(0,0),(1,0),(0,1),(1,1)$.

Their simultaneous kernel is $B_\lambda^j$, because
$\langle Pu,u\rangle=\|a^*u\|^2$, and similarly for $P^c$.
On $\ker P$, the identity $a^*a=I-P$ says that $a$ is an isometry;
its image is $\operatorname{ran}P$, with inverse $a^*$.
Moreover $a$ commutes with $P^c$, again by two anticommutations.
The analogous statements hold for $b$.  It follows that the other three
joint eigenspaces are exactly
$aB_\lambda^{j-1}$, $bB_\lambda^{j-1}$, and
$abB_\lambda^{j-2}$.  Thus
\begin{equation}\label{fockdecomposition}
 E_\lambda^j
 =B_\lambda^j\mathbin{\widehat\oplus}
 aB_\lambda^{j-1}\mathbin{\widehat\oplus}
 bB_\lambda^{j-1}\mathbin{\widehat\oplus}
 abB_\lambda^{j-2},
\end{equation}
where the hat denotes an orthogonal direct sum.  Each indicated map is an
isometry onto its summand; the order in $ab$ affects only a sign.

On these four summands, $\mathcal K$ has respective eigenvalues
$+1,-1,-1,+1$.  If
$\kappa_j(\lambda)=\Tr(\mathcal K_j|_{E_\lambda^j})$, then
\[
 m_j=b_j+2b_{j-1}+b_{j-2},\qquad
 \kappa_j=b_j-2b_{j-1}+b_{j-2}.
\]
For the finite generating polynomials in the degree variable $z$ this is
\[
 m(z)=(1+z)^2b(z),\qquad
 \kappa(z)=(1-z)^2b(z).
\]
Equivalently, as a formal power series at $z=0$,
\begin{equation}\label{fock}
 \kappa(z)=r(z)m(z),\qquad
 r(z)=\frac{(1-z)^2}{(1+z)^2}
      =1+4\sum_{\ell\ge1}(-1)^\ell\ell z^\ell.
\end{equation}
Only finitely many terms of this series enter any fixed degree.

\subsection{Weighted heat coefficients, with all cancellations}

Let $H_j(t)=\Tr(e^{-t\Delta_j})$ and
$W_j(t)=\Tr(\mathcal K_j e^{-t\Delta_j})$.  Multiplying
\eqref{fock} by $e^{-t\lambda}$ and summing over $\lambda>0$ is legitimate
because the heat operators are trace class and $|\kappa_j|\le m_j$.
It gives an exact identity of degree-$p$ heat traces
\begin{equation}\label{exactweightedheat}
 W_p(t)=H_p(t)+4\sum_{\ell=1}^p(-1)^\ell\ell H_{p-\ell}(t)+C_p,
\end{equation}
where $C_p$ is a constant accounting for harmonic forms.
Indeed $\mathcal K$ is the identity on every harmonic space.
In particular, these weighted traces have integer-power heat expansions:
this assertion follows from the exact identity, rather than from a general
claim about arbitrary order-zero pseudodifferential weights.

Write
\[
 H_j(t)\sim(4\pi t)^{-D/2}
 \bigl(a_{0,j}+a_{1,j}t+a_{2,j}t^2+\cdots\bigr).
\]
The ordinary first two coefficients are
\[
 a_{0,j}=\binom Dj V,\qquad
 a_{1,j}=\left(\frac16\binom Dj-\binom{D-2}{j-1}\right)
             \int_Ms\,dV.
\]
Thus their generating polynomials are
\[
 A_0(z)=V(1+z)^D,
\qquad
 A_1(z)=\left(\tfrac16(1+z)^D-z(1+z)^{D-2}\right)\int_Ms\,dV.
\]
Multiplication by $r(z)$ yields
\begin{align*}
 r(z)A_0(z)
 &=V(1-z)^2(1+z)^{D-2}\\
 &=V\bigl((1+z)^D-4z(1+z)^{D-2}\bigr),
\end{align*}
and
\begin{align*}
 r(z)A_1(z)
 &=\left(\tfrac16(1-z)^2(1+z)^{D-2}
          -z(1-z)^2(1+z)^{D-4}\right)\int_Ms\,dV\\
 &=\left(\tfrac16(1+z)^D
          -\tfrac53z(1+z)^{D-2}
          +4z^2(1+z)^{D-4}\right)\int_Ms\,dV.
\end{align*}
Extracting the coefficient of $z^p$, and setting
\[
 N=\binom Dp,\qquad Q=\binom{D-2}{p-1},
 \qquad S_4=\binom{D-4}{p-2},
\]
we obtain
\[
 a_{0,p}^{\mathcal K}=(N-4Q)V,
\qquad
 a_{1,p}^{\mathcal K}
 =\left(\frac N6-\frac{5Q}{3}+4S_4\right)\int_Ms\,dV.
\]
Here and below a binomial coefficient with negative lower index is zero.
The constant $C_p$ in \eqref{exactweightedheat} affects only the
nonnegative powers of $t$; for an exceptional pair $n=D/2\ge3$, it affects
neither of the terms $t^{-n}$ and $t^{1-n}$ used here.

At an exceptional pair, $p(D-p)=D(D-1)/6$, and therefore
\[
 \frac QN=\frac{p(D-p)}{D(D-1)}=\frac16.
\]
Also
\begin{align*}
 \frac{S_4}{N}
 &=\frac{p(p-1)(D-p)(D-p-1)}{D(D-1)(D-2)(D-3)}\\
 &=\frac{p(D-p)\bigl(p(D-p)-(D-1)\bigr)}
        {D(D-1)(D-2)(D-3)}\\
 &=\frac{(D-1)(D-6)}{36(D-2)(D-3)}.
\end{align*}
This remains valid for the smallest pair $(D,p)=(6,1)$, when both sides
vanish.  Substitution gives
\begin{align*}
 \frac1N\left(\frac N6-\frac{5Q}{3}+4S_4\right)
 &=\frac16-\frac5{18}
   +\frac{(D-1)(D-6)}{9(D-2)(D-3)}\\
 &=\frac{-(D-2)(D-3)+(D-1)(D-6)}{9(D-2)(D-3)}\\
 &=-\frac{2D}{9(D-2)(D-3)}.
\end{align*}
We have proved
\begin{equation}\label{weighted}
 a_{0,p}^{\mathcal K}=\frac N3V,\qquad
 a_{1,p}^{\mathcal K}
 =-\frac{2DN}{9(D-2)(D-3)}\int_Ms\,dV.
\end{equation}
The multiplier of total scalar curvature is strictly negative and never
zero at any exceptional pair.

\section{Why the ordinary spectrum determines the weighted coefficient}
\label{sec:recovery}

\subsection{Identifying the otherwise unknown phase}

Equation \eqref{waveK} holds for both $M$ and the model, initially with
possibly different phases $\varepsilon_M$ and $\varepsilon_0$.
Ordinary isospectrality gives, for every $t>0$, the exact equality
\[
 \sum_\lambda m_p(\lambda)e^{i\pi\sqrt\lambda}e^{-t\lambda}
 =\Tr(U_p^M(\pi)e^{-t\Delta_p^M})
 =\Tr(U_p^0(\pi)e^{-t\Delta_p^0}).
\]
These sums converge absolutely.  An order-$-1$ remainder $T$ in
\eqref{waveK} contributes only $o(t^{-D/2})$ to this regularized trace.
Indeed $T(I+\Delta)^{1/2}$ is bounded.  For an orthonormal eigenbasis
$\{\phi_j\}$ with eigenvalues $\lambda_j$,
\[
 |\langle T\phi_j,\phi_j\rangle|
 \le C(1+\lambda_j)^{-1/2}.
\]
Weyl's estimate $\#\{j:\lambda_j\le L\}=O(1+L^{D/2})$ then gives
\[
 |\Tr(Te^{-t\Delta})|
 \le C\sum_j(1+\lambda_j)^{-1/2}e^{-t\lambda_j}
 =O(t^{-(D-1)/2})=o(t^{-D/2})
\]
for $D>1$, as follows also by summing over dyadic spectral intervals.
Hence the leading coefficients in the two traces
are respectively
\[
 (4\pi)^{-D/2}\varepsilon_M a_{0,p}^{\mathcal K}(M),
 \qquad
 (4\pi)^{-D/2}\varepsilon_0 a_{0,p}^{\mathcal K}(\CP^n).
\]
Their $a_0^{\mathcal K}$ values are both $NV/3$ by \eqref{weighted} and
volume equality, and this number is nonzero.  It follows that
\begin{equation}\label{samephase}
 \varepsilon_M=\varepsilon_0=:\varepsilon.
\end{equation}
Notice that this step does not require calculating the Maslov phase
explicitly.

\subsection{A finite-spectrum recovery principle}

The separation argument is not specific to an involution.
We isolate the precise statement, including the data that must be known.

\begin{lemma}[Recovery of a commuting finite-spectrum weight]
\label{lem:finite-weight-recovery}
Let $A$ be a nonnegative self-adjoint elliptic operator of positive order
on a Hermitian vector bundle over a closed manifold.  Let $K$ be a bounded
self-adjoint operator preserving every eigenspace of $A$, and suppose
\[
 \operatorname{Spec}(K)\subseteq S,
\]
where $S\subset\R$ is a fixed finite set.  Let $F$ be a known scalar
function on $\operatorname{Spec}(A)$ and $\varepsilon\ne0$ a known
constant.  Suppose the functional-calculus difference
\[
 T=F(A)-\varepsilon K
\]
extends to a compact operator on $L^2$.  Then, on every sufficiently high
eigenspace $E_\lambda$ of $A$, $K$ is scalar.  Its scalar value is the
unique element $s\in S$ for which $\varepsilon s$ is nearest to
$F(\lambda)$.

Consequently, once $F,\varepsilon,S$ are fixed, the eigenvalues of $A$
with multiplicities determine $\Tr(Ke^{-tA})$ modulo a finite linear
combination of the functions $e^{-t\lambda}$.  In particular, any two
instances with the same ordinary spectral multiset and the same
$F,\varepsilon,S$ have weighted heat traces differing by such a finite
sum.
\end{lemma}

\begin{proof}
Each eigenspace is finite-dimensional.  Since $K$ preserves these
eigenspaces, $T$ does also.  Compactness implies
\begin{equation}\label{eq:finite-weight-tail}
 \|T|_{E_\lambda}\|\longrightarrow0
 \quad\text{as }\lambda\to\infty.
\end{equation}
Indeed, if this failed, there would be pairwise distinct eigenvalues
$\lambda_j\to\infty$ and unit vectors $v_j\in E_{\lambda_j}$ such that
$\|Tv_j\|\ge c>0$.  The vectors $v_j$ are orthonormal, hence weakly
converge to zero; a compact operator sends such a sequence to a
norm-null sequence.  This contradicts the lower bound.

If $S$ is a singleton, the conclusion is immediate.  Otherwise put
\[
 d=\min\{|\varepsilon r-\varepsilon s|:r,s\in S,\ r\ne s\}>0.
\]
Choose a threshold above which the norms in
\eqref{eq:finite-weight-tail} are less than $d/3$.  If $s$ is an
eigenvalue of the self-adjoint restriction $K|_{E_\lambda}$, then
\[
 |F(\lambda)-\varepsilon s|<d/3.
\]
There cannot be two distinct such elements, because their centers are
at least distance $d$ apart.  Since the restriction is self-adjoint and
has only one eigenvalue, it equals $sI$.  Every other center has distance
greater than $2d/3$ from $F(\lambda)$, so $s$ is uniquely nearest.

Only finitely many eigenvalues lie below the threshold, each with finite
multiplicity.  All high weighted multiplicities are therefore the
known numbers $s(\lambda)\dim E_\lambda$, and any two possible weighted
traces differ only in the finitely many low eigenspaces.  The heat
operators are trace class because $A$ is elliptic of positive order on a
closed manifold, and boundedness of $K$ ensures that the weighted traces
are defined for every $t>0$.
\end{proof}

For an involution, $S=\{-1,1\}$.  In the application here
$F(\lambda)=e^{i\pi\sqrt\lambda}$ and $|\varepsilon|=1$, so the
recovered high-energy sign is explicitly
\begin{equation}\label{eq:explicit-high-sign}
 s(\lambda)=\operatorname{sgn}
       \operatorname{Re}(\varepsilon^{-1}e^{i\pi\sqrt\lambda}).
\end{equation}
The real part is bounded away from zero at sufficiently high energy by
the compact-remainder hypothesis.  A finite sum of exponentials is
analytic at $t=0$, so all negative-power terms of a weighted heat
expansion are determined by the ordinary spectrum.  More generally,
any logarithmic or other nonanalytic coefficients in a specified
weighted expansion are unchanged by the finite ambiguity.

\begin{lemma}[Stability of recovered signs under matched phases]
\label{lem:finite-weight-matched-phases}
Consider two instances of Lemma~\ref{lem:finite-weight-recovery} with
the same finite set $S$ and the same $\varepsilon$.  Suppose their
ordinary eigenvalues can be paired, with multiplicities, as
$\lambda_j$ and $\mu_j$ tending to infinity and satisfying
\[
 F(\lambda_j)-F(\mu_j)\longrightarrow0.
\]
Then their recovered weight values agree on all sufficiently high
paired indices.  This conclusion does not require that a multiple
eigenvalue on one side remain a single eigenvalue on the other.
\end{lemma}
\begin{proof}
For $|S|\ge2$, let $d$ be the minimum center separation in the preceding
proof.  Choose a tail on which both compact error bounds and
$|F(\lambda_j)-F(\mu_j)|$ are less than $d/4$.  Distinct recovered
values $r,s$ would imply
\[
 |\varepsilon r-\varepsilon s|
 \le |\varepsilon r-F(\lambda_j)|
       +|F(\lambda_j)-F(\mu_j)|
       +|F(\mu_j)-\varepsilon s|<3d/4,
\]
contradicting its definition.  The singleton case is immediate.
\end{proof}

\subsection{Application to the K\"ahler involution}

Apply Lemma~\ref{lem:finite-weight-recovery} with
$A=\Delta_p$, $K=\mathcal K_p$, $S=\{-1,1\}$, and
$F(\lambda)=e^{i\pi\sqrt\lambda}$. Equation~\eqref{waveK} gives
the compact remainder, and \eqref{samephase} gives the same
$\varepsilon$ on the two isospectral manifolds. The lemma therefore
shows that the high weighted multiplicities agree exactly, even
though only the ordinary spectrum in the one fixed degree is given.
Consequently, for some finite threshold $L$,
\begin{equation}\label{finiteweighted}
 W_p^M(t)-W_p^0(t)=\sum_{\lambda\leq L}c_\lambda e^{-t\lambda}
\end{equation}
with finitely many integer coefficients $c_\lambda$.

The right side is analytic at $t=0$. Since an exceptional pair has
$n\geq3$, it cannot alter the negative-power term $t^{1-n}$.
The exact heat identity \eqref{exactweightedheat} gives the required
expansions, so \eqref{weighted} implies
\[
\int_Ms\,dV=\int_{\CP^n}s_0\,dV.
\]
Thus \eqref{meanscalar} gives $k=n+1$, and the
Kobayashi--Ochiai argument identifies the complex manifold and class
also in every exceptional case. The lower-degree heat traces in
\eqref{exactweightedheat} were used only to calculate universal
coefficients; no isospectrality in those additional degrees was assumed.

\section{The metric is determined after the complex manifold and class}
\label{sec:metric}

By the preceding argument, there is a biholomorphism between $M$ and
$\CP^n$, and under this identification $[\omega]=[\omega_0]=\pi h$.
In this section we compare the two metrics on this fixed complex manifold.
Write $dV=\omega^n/n!$, and use $\delta I=I(g)-I(g_0)$ for the
difference of a curvature integral. The measures in the two terms are
their respective Riemannian volume measures.

\subsection{Curvature conventions and the Chern--Weil identities}

Let $R_c$ be the complex K\"ahler curvature tensor, and write
$\Ric_c$ for its complex Ricci contraction. We use the norms satisfying
\begin{equation}\label{realcomplexnorms}
 |\operatorname{Rm}|^2=4|R_c|^2,\qquad
 |\Ric_g|^2=2|\Ric_c|^2,\qquad
 \operatorname{tr}_{\C}\Ric_c=s/2.
\end{equation}
These are the conventions in \cite[Lemmas 3.3--3.5]{Li} and
\cite[Section 2]{HLXZ}. Our K\"ahler form is the ordinary metric form
$\omega(X,Y)=g(JX,Y)$; hence the Ricci form $\rho$ represents
$2\pi c_1(M)$.

For $n\geq2$, the Apte identities, in this normalization, are
\begin{align}
 \int_M\left(\frac{s^2}{4}-|\Ric_c|^2\right)dV
 &=\frac{4\pi^2}{(n-2)!}
   \left\langle c_1(M)^2[\omega]^{n-2},[M]\right\rangle,
   \label{apte1}\\
 \int_M\left(\frac{s^2}{4}-2|\Ric_c|^2+|R_c|^2\right)dV
 &=\frac{8\pi^2}{(n-2)!}
   \left\langle c_2(M)[\omega]^{n-2},[M]\right\rangle.
   \label{apte2}
\end{align}
We recall their algebraic origin to make the normalizations explicit.
For real $(1,1)$-forms $\alpha,\beta$, diagonalization in a unitary frame
gives the wedge-contraction identity
\[
 \frac{\alpha\wedge\beta\wedge\omega^{n-2}}{(n-2)!}
 =\bigl((\Lambda\alpha)(\Lambda\beta)
               -\langle\alpha,\beta\rangle\bigr)dV.
\]
Taking $\alpha=\beta=\rho$, using
$\Lambda\rho=s/2$ and $|\rho|^2=|\Ric_c|^2$, proves
\eqref{apte1}. For the second identity, let $\Theta$ be the curvature
matrix of the Chern connection of $T^{1,0}M$, and put $A=i\Theta$.
The Chern--Weil representatives are
\[
 c_1=\frac{\operatorname{Tr}A}{2\pi},\qquad
 c_2=\frac{(\operatorname{Tr}A)\wedge(\operatorname{Tr}A)
                   -\operatorname{Tr}(A\wedge A)}{8\pi^2}.
\]
Applying the wedge-contraction identity entrywise gives
\[
 \frac{\operatorname{Tr}(A\wedge A)\wedge\omega^{n-2}}{(n-2)!}
       =\bigl(|\Ric_c|^2-|R_c|^2\bigr)dV;
\]
substitution gives \eqref{apte2}. See \cite{Apte,Li} for the equivalent
classical formulations.

The right sides of \eqref{apte1}--\eqref{apte2} agree for the two
metrics, because both the complex manifold and the K\"ahler class
agree. Subtracting first \eqref{apte1} and then \eqref{apte2} yields
\begin{equation}\label{fixedclassdifferences}
 \delta\!\int_M|\Ric_c|^2dV
 =\delta\!\int_M|R_c|^2dV
 =\frac14\delta\!\int_Ms^2dV.
\end{equation}
When $n=1$, no $c_2$ identity is required: the equalities
$|R_c|^2=|\Ric_c|^2=s^2/4$ hold pointwise, so
\eqref{fixedclassdifferences} remains valid.

\subsection{The second heat coefficient forces constant scalar curvature}

Put
\[
 Q=\binom{D-2}{p-1},\qquad S_4=\binom{D-4}{p-2},
\]
with $Q=0$ if $p=0$ and $S_4=0$ if $p<2$.
The standard integrated Hodge heat coefficient is
\begin{equation}\label{a2full}
 a_{2,p}=\int_M
    \bigl(4\lambda_1|R_c|^2+2\lambda_2|\Ric_c|^2
                   +\lambda_3s^2\bigr)dV,
\end{equation}
where \cite{Patodi,Gilkey,Li}
\begin{align*}
 \lambda_1&=\frac{N}{180}-\frac{Q}{12}+\frac{S_4}{2},\\
 \lambda_2&=-\frac{N}{180}+\frac{Q}{2}-2S_4,\\
 \lambda_3&=\frac{N}{72}-\frac{Q}{6}+\frac{S_4}{2}.
\end{align*}
The integral of the scalar Laplacian term in the local coefficient is
zero because $M$ is closed. In the combination needed here the
$Q$ and $S_4$ terms cancel:
\begin{align*}
 \lambda_1+\frac{\lambda_2}{2}+\lambda_3
 &=N\left(\frac1{180}-\frac1{360}+\frac1{72}\right)
   +Q\left(-\frac1{12}+\frac14-\frac16\right)\\
 &\hspace{15mm}+S_4\left(\frac12-1+\frac12\right)
 =\frac{N}{60}>0.
\end{align*}
The calculation is valid beyond projective space and gives the
following general identity.

\begin{proposition}[The second heat coefficient in a fixed K\"ahler class]
\label{prop:fixed-class-a2}
Let $g,g_*$ be K\"ahler metrics on the same compact complex manifold
$X$ of complex dimension $n\geq1$, with $[\omega_g]=[\omega_{g_*}]$.
For every $0\leq p\leq2n$,
\begin{equation}\label{eq:fixed-class-calabi}
a_{2,p}(g)-a_{2,p}(g_*)
=\frac{\binom{2n}{p}}{60}
\left(\int_Xs_g^2dV_g-\int_Xs_{g_*}^2dV_{g_*}\right).
\end{equation}
If $g_*$ has constant scalar curvature $s_*$, then
\begin{equation}\label{eq:fixed-class-variance}
a_{2,p}(g)-a_{2,p}(g_*)
=\frac{\binom{2n}{p}}{60}\int_X(s_g-s_*)^2dV_g.
\end{equation}
\end{proposition}
\begin{proof}
For $n\geq2$, the two Chern--Weil pairings in
\eqref{apte1}--\eqref{apte2} depend only on $X$ and the common class.
Thus \eqref{fixedclassdifferences} holds with $g_*$ in place of $g_0$.
Substitution into \eqref{a2full} and the cancellation just proved yield
\eqref{eq:fixed-class-calabi}. This works directly for $p\leq n$;
Hodge duality gives the remaining degrees. When $n=1$,
$a_{2,0}=\int_Xs^2dV/60$ and $a_{2,2}=a_{2,0}$.
The identity $H_1(t)=2H_0(t)-\chi(X)$ gives $a_{2,1}=2a_{2,0}$,
because its constant Euler-characteristic term does not affect the
$t^1$ coefficient. This proves the formula also for surfaces.

The common K\"ahler class fixes volume and total scalar curvature:
\[
V=\frac{\langle[\omega]^n,[X]\rangle}{n!},\qquad
\int_Xs\,dV=\frac{4\pi}{(n-1)!}
\langle c_1(X)[\omega]^{n-1},[X]\rangle.
\]
If $s_*$ is constant, the latter integral equals $s_*V$, and expanding
$\int_X(s_g-s_*)^2dV_g$ proves \eqref{eq:fixed-class-variance}.
\end{proof}

We now apply this identity to the projective-space metrics already
identified above.
Combining isospectral equality of $a_{2,p}$ with
\eqref{fixedclassdifferences} and \eqref{a2full} therefore gives
\[
 0=\delta a_{2,p}=\frac{N}{60}\delta\!\int_Ms^2dV.
\]
We already know $\Vol(g)=\Vol(g_0)=V$ and $\int_Ms\,dV=s_0V$.
Consequently
\[
 \delta\!\int_Ms^2dV
 =\int_M(s-s_0)^2dV=0.
\]
The integrand is continuous and nonnegative, so $s\equiv s_0$.
Notice that this argument makes no sign assumption on any individual
$\lambda_i$; it also applies to $(p,n)=(2,8)$, where $\lambda_1=0$.

\subsection{The traceless Ricci and Bochner tensors vanish}

Let $\Ric_c^0=\Ric_c-(s/(2n))g_c$. The orthogonal decomposition gives
\[
 |\Ric_c|^2=|\Ric_c^0|^2+\frac{s^2}{4n}.
\]
Since \eqref{fixedclassdifferences} and $s=s_0$ imply equality of the
integrated Ricci norms, and the model has $\Ric_c^0=0$, we obtain
\[
 \int_M|\Ric_c^0|^2dV
 =\int_M|\Ric_c|^2dV-\frac{s_0^2}{4n}V=0.
\]
Thus the metric is K\"ahler--Einstein. This also supplies a direct
proof of the usual implication from constant scalar curvature to
Einstein in the present cohomological situation.

For $n\geq2$, write the orthogonal K\"ahler curvature decomposition as
$R_c=\mathcal S+\mathcal P+\mathcal B$, where $\mathcal B$ is the
Bochner tensor. The explicit formulas are
\begin{align*}
 \mathcal S_{i\bar j k\bar l}
 &=\frac{s}{2n(n+1)}
     (g_{i\bar j}g_{k\bar l}+g_{i\bar l}g_{k\bar j}),\\
 \mathcal P_{i\bar j k\bar l}
 &=\frac1{n+2}
    (g_{i\bar j}(\Ric_c^0)_{k\bar l}
     +g_{k\bar l}(\Ric_c^0)_{i\bar j}
     +g_{i\bar l}(\Ric_c^0)_{k\bar j}
     +g_{k\bar j}(\Ric_c^0)_{i\bar l}),
\end{align*}
and their contractions give
\[
 |R_c|^2=\frac{s^2}{2n(n+1)}
            +\frac4{n+2}|\Ric_c^0|^2+|\mathcal B|^2.
\]
The model has $\mathcal P=\mathcal B=0$. Equality of integrated
$|R_c|^2$ in \eqref{fixedclassdifferences}, together with $s=s_0$ and
$\Ric_c^0=0$, now forces $\int_M|\mathcal B|^2dV=0$. Hence
$R_c=\mathcal S$ everywhere, which says that $g$ has constant
holomorphic sectional curvature $4$. For $n=1$, $s=s_0=8$ already says
that the Gaussian, equivalently holomorphic sectional, curvature is $4$.

Finally, a complete simply connected K\"ahler manifold of positive
constant holomorphic sectional curvature is holomorphically isometric
to the corresponding Fubini--Study space. This is the uniformization theorem. Compactness gives
completeness, and \eqref{topology} gives simple connectedness. Thus
Theorem~\ref{main} follows.


\section{Consequences: asymptotic spectra, spectral tails, and Besse metrics}
\label{sec:consequences}

We retain the normalization in which the model has holomorphic sectional
curvature $4$, K\"ahler class $\pi h$, volume
$V_0=\pi^n/n!$, and scalar curvature $s_0=4n(n+1)$.
All spectra below include multiplicities.  The conclusions for any
other Fubini--Study scale follow by simultaneously rescaling the two
metrics and their eigenvalues.

The statements are given for every form degree $0\le p\le2n$.
Hodge duality identifies the degree-$p$ spectrum with the degree-$(2n-p)$
spectrum, so all arguments reduce to
$p'=\min\{p,2n-p\}\le n$.  The quantities
\[
 N=\binom{2n}{p},\qquad
 q=p^2-2np+\frac{n(2n-1)}3
\]
are invariant under this replacement.  Thus the distinction between
ordinary and exceptional degrees is unchanged.

\subsection{Rigidity under asymptotic matching of eigenvalue lists}

\begin{theorem}[Asymptotic spectral rigidity]
\label{thm:asymptotic-rigidity}
Let $n\ge2$ and $0\le p\le2n$.  Let $(M,g,J)$ be a compact connected
K\"ahler $n$-fold without boundary.  List the degree-$p$ Hodge eigenvalues
of $M$ and normalized Fubini--Study space in nondecreasing order,
including multiplicities, as $\lambda_j$ and $\mu_j$, with the same
indexing convention.  If
\begin{equation}\label{eq:asymptotic-matching}
 |\lambda_j-\mu_j|=o(\mu_j^{-1})\qquad(j\longrightarrow\infty),
\end{equation}
then $(M,g,J)$ is holomorphically isometric to the model.
\end{theorem}

\begin{proof}
By duality we may work in degree $p'\le n$ and write it as $p$.
Set $\delta_j=\lambda_j-\mu_j$.  The expression
\eqref{eq:asymptotic-matching} concerns only the positive high-energy
tail; the finitely many zero eigenvalues cause no ambiguity.

\emph{Step 1: recovering the first three ordinary heat coefficients.}
The model's counting function satisfies
$N_0(R)=O(1+R^n)$ by Weyl's law.  For any fixed $c>0$ and
$0<a<n$, summation over unit or dyadic spectral intervals gives
\begin{equation}\label{eq:weighted-Weyl-sum}
 \sum_{\mu_j\ge1}\mu_j^{-a}e^{-ct\mu_j}
 =O(t^{a-n})\qquad(0<t\le1).
\end{equation}
For example, this follows by comparing the sum with a
Stieltjes integral against $N_0(R)$ and integrating by parts, using
$N_0(R)=O(R^n)$ for $R\ge1$.

For sufficiently large $j$, the two eigenvalues are comparable, say
$\lambda_j\ge\mu_j/2$.  The mean-value theorem then gives
\[
 |e^{-t\lambda_j}-e^{-t\mu_j}|
 \le t|\delta_j|e^{-ct\mu_j}
\]
with a fixed $c>0$.  Given $\eta>0$, choose a tail on which
$|\delta_j|\le\eta\mu_j^{-1}$.  Since $n>1$, formula
\eqref{eq:weighted-Weyl-sum} with $a=1$ bounds the sum of this tail
by $C\eta t^{2-n}$.  For each of the finitely many earlier paired
indices, $e^{-t\lambda_j}-e^{-t\mu_j}=O(t)$ as $t\to0$; the constant
terms cancel because the indices are paired one for one.  Since
$O(t)=o(t^{2-n})$ for $n\ge2$, letting $\eta\to0$ proves
\begin{equation}\label{eq:asymptotic-heat-small}
 H_{p,M}(t)-H_{p,0}(t)=o(t^{2-n}).
\end{equation}
Comparison with the heat expansions implies
\begin{equation}\label{eq:asymptotic-a012}
 a_{r,p}(M,g)=a_{r,p}(\CP^n,g_0),\qquad r=0,1,2.
\end{equation}
In particular this includes the constant heat-trace coefficient
$a_{2,p}$ when $n=2$.

\emph{Step 2: periodicity and projective cohomology.}
For large $j$,
\[
 |\sqrt{\lambda_j}-\sqrt{\mu_j}|
 =\frac{|\delta_j|}{\sqrt{\lambda_j}+\sqrt{\mu_j}}
 =o(\mu_j^{-3/2}).
\]
The model spectral branches satisfy
$e^{2\pi i\sqrt{\mu_j}}\to1$.  Hence the same limit holds with
$\lambda_j$ in place of $\mu_j$, and
$U_M(2\pi)-I$ is compact by the spectral compactness criterion.
The Egorov argument gives \eqref{period2} for $M$.  The universal-cover
and Besse-cohomology arguments leading to \eqref{topology} use only
this common period and the K\"ahler condition, so they apply unchanged.
Volume equality in \eqref{eq:asymptotic-a012} then fixes
$[\omega]=\pi x$ and gives the pairing \eqref{meanscalar}.

If $q\ne0$, equality of $a_{1,p}$ already yields
$\int_Ms\,dV=s_0V_0$.  The maximal-index argument identifies the
complex manifold with $\CP^n$, and equality of $a_{2,p}$ completes
the metric argument in Section~\ref{sec:metric}.  It remains to
consider $q=0$.

\emph{Step 3: the leading return traces still agree in the exceptional case.}
Every exceptional pair has $n\ge3$.  Set
\[
 F_t(x)=e^{i\pi\sqrt{x}}e^{-tx}\qquad(x>0).
\]
Its derivative satisfies
\[
 |F_t'(x)|\le
 \left(\frac{\pi}{2\sqrt{x}}+t\right)e^{-tx}.
\]
On the high tail, comparability of $\lambda_j$ and $\mu_j$ and the
mean-value theorem therefore yield
\[
 |F_t(\lambda_j)-F_t(\mu_j)|
 \le C|\delta_j|(\mu_j^{-1/2}+t)e^{-ct\mu_j}.
\]
The hypothesis implies the weaker uniform bound
$|\delta_j|\le C\mu_j^{-1}$ on a tail.  Using
\eqref{eq:weighted-Weyl-sum} with $a=3/2$ and $a=1$ gives the explicit
rough estimate
\begin{align}
 &\Tr(U_{p,M}(\pi)e^{-t\Delta_{p,M}})
   -\Tr(U_{p,0}(\pi)e^{-t\Delta_{p,0}})\notag\\
 &\hspace{15mm}
 =O\bigl(1+t^{3/2-n}+t^{2-n}\bigr)
 =o(t^{-n}).\label{eq:asymptotic-wave-leading}
\end{align}
The finitely many low-energy terms contribute $O(1)$; no derivative
estimate at $x=0$ is needed.

The proof of the exceptional return and holonomy results uses ordinary
isospectrality only to compare the leading $t^{-n}$ coefficients of
these traces, in addition to the compact period already established.
Equation \eqref{eq:asymptotic-wave-leading} supplies precisely that
comparison.  More explicitly, since $G_M^{2\pi}=I$, $G_M^\pi$ is an
involution with clean fixed set.  The nonzero model leading coefficient
rules out a proper fixed set, giving \eqref{periodpi}.  Compactness of
$U_M(2\pi)-I$ then forces involutive return holonomy.  Equality of the
leading coefficients and equality of volumes give the same exterior
character equation \eqref{character}.  The Krawtchouk estimate and
projective-bundle obstruction depend only on this equation and
\eqref{topology}; hence they again yield \eqref{modelhol}.
Consequently \eqref{waveK} holds for both manifolds.

The leading term in \eqref{waveK} is
$\varepsilon a_{0,p}^{\mathcal K}$, and
\eqref{weighted} gives $a_{0,p}^{\mathcal K}=NV_0/3\ne0$ on both sides.
Thus \eqref{eq:asymptotic-wave-leading} also identifies the same phase
$\varepsilon$ in the two compact-remainder relations.  Equality of
lower ordinary wave coefficients is not asserted or needed.

\emph{Step 4: recovering the weighted first heat coefficient.}
Apply Lemma~\ref{lem:finite-weight-recovery} to the exact involutions
$\mathcal K_p$ and the function $F(\lambda)=e^{i\pi\sqrt\lambda}$.
Each sufficiently high eigenspace has a single $\mathcal K_p$-sign.
Moreover,
\[
 |F(\lambda_j)-F(\mu_j)|
 \le\pi|\sqrt{\lambda_j}-\sqrt{\mu_j}|\longrightarrow0.
\]
Lemma~\ref{lem:finite-weight-matched-phases} therefore shows that the
two signs agree on all sufficiently high paired indices.  This argument
allows one model eigenspace to split into several nearby eigenvalues on
$M$: it uses each repeated index and its phase, not an assumed
identification of whole eigenspaces.

Write the common high sign as $\sigma_j\in\{-1,1\}$.  For a fixed
large $J$,
\[
 W_{p,M}(t)-W_{p,0}(t)
 =\sum_{j>J}\sigma_j
       (e^{-t\lambda_j}-e^{-t\mu_j})+O(1).
\]
The absolute-value estimate in Step~1 bounds the tail by
$o(t^{2-n})$.  Since $n\ge3$ in the exceptional case,
\begin{equation}\label{eq:asymptotic-weighted-small}
 W_{p,M}(t)-W_{p,0}(t)
 =o(t^{2-n})+O(1)=o(t^{1-n}).
\end{equation}
Comparison of the weighted expansions, whose leading coefficients
already agree, gives equality of $a_{1,p}^{\mathcal K}$.
Its multiplier of total scalar curvature in \eqref{weighted} is
nonzero.  Thus $\int_Ms\,dV=s_0V_0$ also when $q=0$.

Equation \eqref{meanscalar} now gives $k=n+1$.  The
Kobayashi--Ochiai theorem identifies the complex manifold, and the
fixed-class calculation of Section~\ref{sec:metric}, using the ordinary
$a_{2,p}$ equality from Step~1, identifies the metric.  This completes
the proof.
\end{proof}

The little-oh assumption in \eqref{eq:asymptotic-matching} is used to
eliminate a possible defect at the $a_2$ scale.  Replacing it by
$O(\mu_j^{-1})$ yields only $O(t^{2-n})$ in
\eqref{eq:asymptotic-heat-small} and does not prove equality of
$a_{2,p}$.  In complex dimension one, the relevant Weyl sum with
$a=1$ is logarithmic, so the preceding proof does not establish that
endpoint.  These remarks describe limitations of this proof, rather
than counterexamples to weaker hypotheses.

\subsection{Rigidity from an exact spectral tail}

\begin{proposition}[Spectral-tail rigidity]
\label{prop:spectral-tail-rigidity}
Let $n\ge3$ and $0\le p\le2n$.  Suppose the degree-$p$ Hodge spectral
multisets of a compact connected K\"ahler $n$-fold and normalized
Fubini--Study space agree above a finite threshold.  Then the two
manifolds are holomorphically isometric.
\end{proposition}

\begin{proof}
After Hodge duality, work in degree at most $n$.  The ordinary heat
traces differ by a finite signed sum
\begin{equation}\label{eq:finite-tail-defect}
 H_{p,M}(t)-H_{p,0}(t)
 =\sum_{\lambda}c_\lambda e^{-t\lambda},
\end{equation}
where $c_\lambda$ is the multiplicity on $M$ minus that on the model,
and only finitely many coefficients are nonzero.  This is analytic
at zero.  For $n\ge3$, the terms involving $a_{0,p},a_{1,p},a_{2,p}$
have powers $t^{-n},t^{1-n},t^{2-n}$, all negative.  Thus these three
coefficients agree.

The compact-period property depends only on the high spectrum and is
unchanged by finitely many eigenvalues.  The periodicity and topology
arguments therefore give \eqref{period2}, \eqref{topology}, and the
same class normalization.  If $q\ne0$, equality of $a_{1,p}$ and
$a_{2,p}$ now completes the argument as before.

If $q=0$, the regularized $\pi$-wave traces also differ by a finite
sum, namely
\[
 \sum_{\lambda}c_\lambda e^{i\pi\sqrt\lambda}e^{-t\lambda}=O(1).
\]
Their leading $t^{-n}$ coefficients agree.  This gives the shorter
return, the character equation, and the model holonomy by exactly
the exceptional argument used in Theorem~\ref{main}.  It also fixes
the common phase in \eqref{waveK}.  The finite-spectrum recovery lemma
then identifies the $\mathcal K_p$-signs and weighted multiplicities
at every sufficiently high common eigenvalue.  The weighted heat
traces consequently differ by another finite exponential sum.
Their $t^{1-n}$ coefficients agree, and \eqref{weighted} recovers
total scalar curvature.  The maximal-index and final metric
arguments again apply.
\end{proof}

Exact equality of spectral tails is stated separately from
Theorem~\ref{thm:asymptotic-rigidity}.  Different numbers of low
eigenvalues may shift the two ordered lists by a nonzero number of
indices.  Thus tail equality does not in general provide the
same-index matching in \eqref{eq:asymptotic-matching}.

There are also sufficient low-dimensional versions.  Using the signed
defect coefficients in \eqref{eq:finite-tail-defect}, put
\[
 C_0=\sum_\lambda c_\lambda,
 \qquad C_1=\sum_\lambda c_\lambda\lambda.
\]
Its Taylor expansion is $C_0-C_1t+O(t^2)$.  In complex dimension two,
the same rigidity conclusion follows from tail equality plus
$C_0=0$: the negative-power coefficients $a_0,a_1$ are automatic,
and $C_0=0$ preserves the constant coefficient $a_2$.
In complex dimension one, the conditions $C_0=C_1=0$ preserve the
constant and linear coefficients, namely $a_1,a_2$; $a_0$ is
automatic.  The compact-period and topology arguments still apply.
There are no exceptional pairs in dimensions one or two, so the
ordinary first coefficient recovers total scalar curvature and the
metric argument finishes these cases.  These moment conditions are
sufficient; no necessity claim is made.

\subsection{A criterion using Besse geometry and three numerical invariants}

\begin{theorem}[Besse and heat-invariant criterion]
\label{thm:besse-heat-criterion}
Let $n\ge1$, $0\le p\le2n$, and let $(M,g,J)$ be a compact connected
K\"ahler $n$-fold without boundary.  Suppose its unit-speed geodesic
flow has a common period and
\begin{equation}\label{eq:besse-numerical-data}
 \Vol(M,g)=V_0,\qquad
 \int_Ms\,dV=s_0V_0,\qquad
 a_{2,p}(M,g)=a_{2,p}(\CP^n,g_0).
\end{equation}
Then $(M,g,J)$ is holomorphically isometric to normalized
Fubini--Study space.
\end{theorem}

\begin{proof}
A common geodesic period, together with the K\"ahler condition,
supplies all hypotheses of the topology argument leading to
\eqref{topology}.  In particular, $M$ is simply connected and
$H^*(M;\Z)=\Z[x]/(x^{n+1})$.  Choose the positive generator as
before.  Volume equality determines $[\omega]=\pi x$, and write
$c_1(M)=kx$.  The pairing \eqref{meanscalar} and the prescribed
total scalar curvature imply
\[
 4nk=s_0=4n(n+1),\qquad k=n+1.
\]
The positive integral class $x$ is the first Chern class of an ample
line bundle, so the maximal-index theorem identifies $M$ with
$\CP^n$.  The normalized K\"ahler classes agree.  In that fixed
class, the last equality in \eqref{eq:besse-numerical-data} forces
scalar-curvature variance, traceless Ricci curvature, and Bochner
curvature to vanish by Section~\ref{sec:metric}.  Hence the resulting
holomorphic isometry has precisely the stated scale.
Hodge duality treats degrees greater than $n$ without changing any
of these numerical assumptions.
\end{proof}

For $q\ne0$, equality with the model of $a_{0,p},a_{1,p},a_{2,p}$
implies all three conditions in \eqref{eq:besse-numerical-data}.
Thus in a Besse K\"ahler class the three ordinary coefficients are
sufficient in the nonexceptional degrees.  At $q=0$, the theorem
prescribes total scalar curvature explicitly; it does not claim to
recover that quantity from the three ordinary coefficients alone.

The common-period hypothesis may be replaced, for the
purpose of this criterion, by the following sufficient spectral
condition: for some $T>0$ and $|\zeta|=1$,
\begin{equation}\label{eq:phase-clustering-criterion}
 e^{iT\sqrt\lambda}\longrightarrow\zeta
 \quad\text{as }\lambda\to\infty
 \text{ through }\Spec_p(M,g).
\end{equation}
Indeed, the diagonal multipliers of
$U_p(T)-\zeta I$ tend to zero, so this operator is compact.
For every order-zero pseudodifferential operator $A$ with scalar
principal symbol, the scalar phase cancels in
$U_p(T)^*AU_p(T)-A$, which is therefore compact.  Egorov's theorem
then gives $a\circ G^T=a$ for all scalar symbols $a$, forcing
$G^T=I$.  This recovers the geometric hypothesis of
Theorem~\ref{thm:besse-heat-criterion}.

The phase condition \eqref{eq:phase-clustering-criterion} alone does
not supply the shorter-return character comparison used in the
exceptional proof.  Thus neither that condition nor the three
ordinary heat coefficients are being substituted for the explicit
total-scalar-curvature hypothesis at $q=0$.

\end{document}